\documentclass[11pt]{article}
\usepackage{amssymb,amsmath}   

\usepackage{graphicx}

\def\mb#1{{\mbox{\boldmath $#1$}}}

\newcommand{\be}{\begin{eqnarray}}
\newcommand{\nd}{\end{eqnarray}}

\begin{document}
\begin{center}
{\bf On uniform continuity of Cauchy's function and\\ 
	uniform convergence of Cauchy's integral formula with applications}\\
 by Theodore Yaotsu Wu\footnote{California Institute of Technology, Pasadena, CA 91125 U.S.A. ~Email: tywu@caltech.edu}
\end{center}
\vskip 1.8mm
\noindent{\bf Abstract.}  ~This study is on Cauchy's function $f(z)$ and its integral, $(2\pi i)J[f(z)]\equiv \oint_C f(t)dt/(t-z)$ taken along a closed simple contour $C$, in regard to their comprehensive properties over the entire $z=x+iy$ plane consisted of the open domain ${\cal D}^+$ bounded by $C$ and the open domain ${\cal D}^-$ outside $C$.  (i) With $f(z)$ assumed to be $C^n$ ($n$-times continuously differentiable) $\forall z\in {\cal D}^+$ and in a neighborhood of $C$, $f(z)$ and its derivatives $f^{(n)}(z)$ are proved uniformly continuous in the closed domain $\overline{{\cal D}^+}=[\cal D^++C]$.  (ii) Under this new assumption, Cauchy's integral formula and its derivatives $\forall z\in {\cal D}^+$ (or $\forall z\in {\cal D}^-$) are proved to converge uniformly in $\overline{{\cal D}^+}$ (or in $\overline{{\cal D}^-} =[{\cal D}^-+C]$), respectively, thereby rendering the integral formula valid over the entire $z$-plane.  (iii) The same claims (as for $f(z)$ and $J[f(z)]$) are shown extended to hold for the complement function $F(z)$, defined to be $C^n \forall z\in {\cal D}^-$ and about $C$.  (iv) The uniform convergence theorems for $f(z)$ and $F(z)$ shown for arbitrary contour $C$ are adapted to apply to special domains in the upper or lower half $z$-planes and those inside and outside the unit circle $|z|=1$ to achieve the generalized Hilbert transforms for these cases.  (v) Further, the singularity distribution of $f(z)$ in ${\cal D}^-$ is elucidated by considering the direct problem exemplified with several typical singularities prescribed in ${\cal D}^-$.  (vi) A comparative study is made between generalized integral formulas and Plemelj's formulas on their basic properties.  (vii) Physical significances of these formulas are illustrated with applications to nonlinear airfoil theory.  (viii) Finally, an unsolved inverse problem to determine all the singularities of Cauchy function $f(z)$ in domain ${\cal D}^-$ is presented for resolution as a conjecture.
\vskip 1mm
\noindent{\bf Key word:} ~Cauchy function, Cauchy integral, uniform continuity, uniform convergence, integral transforms, functional properties. 
\vskip 2mm
\noindent{\bf 1. ~Introduction.}
\vskip 0.8mm
In the classical studies of water waves propagating in permanent form on water of finite or infinite depth in the physical $z=x+iy$ plane, the complex potential $f=\phi+i\psi$ was first adopted by Stokes (1880)[5] as the independent variable for the advantage due to its boundary location being known.  The infinite strip of $f$ for a solitary wave or a semi-infinite strip of $f$ for a periodic wave train on water of infinite depth can readily be mapped conformally onto a unit disc ($|\zeta|\leq 1$) in a $\zeta=\xi+i\eta$ plane.  For waves of a rounded crest, the complex velocity $w(z(\zeta))=df/dz$ is an analytic function, regular (or holomorphic) in the closed flow domain.  As the solitary, or a periodic wave grows to become the highest with a corner crest of interior angle of $120^\circ$, $w(z(\zeta))$ then has a primary (Stokes's) and secondary (Grant's) algebraic branch singularities at the crest on the flow boundary $|\zeta|=1$, whereas $w(z(\zeta))$ remains regular in the open domain $|\zeta|<1$ of the inner flow field. 

The first and foremost query is to question: ``Where do the branch singularities go as the highest wave becomes rounded at the crest by a slight reduction in height?"  ``Do the singularities retain their type, with only changes in strength and location away from the flow boundary to an interior point in the open domain ${\cal D}^-\{\forall \zeta \in |\zeta|>1\} $ outside the disc? or even with changes in types of singularities?"  In fact, for waves that are regular in the closed flow domain $\overline{{\cal D}^+}\{\forall \zeta \in |\zeta|\leq 1\}$, $w(z(\zeta))$ must have singularities distributed in the open domain ${\cal D}^-$, for otherwise, by Liouville's Theorem, $w(z(\zeta))$ would have to be a constant in the entire $\zeta$-plane, which is the trivial case of a wave of vanishing amplitude.

To determine the singularities of $w(z(\zeta))$ in domain ${\cal D}^-$ outside the flow field, to be called {\it the inverse problem} in short, is thought of as essential to gaining sound physical explanations for a broad list of curious properties of these waves, known or newly discovered.  This acute need has been greatly motivated by a recent series of our studies, in which Part 1[8] has developed a {\it unified intrinsic functional expansion theory} for exact evaluation on the Euler model for solitary waves of arbitrary height $a$ on water of depth $h$.  This theory adopts an expansion in terms of a set of intrinsic component functions determined from analysis of the flow velocity field about the wave crest, its outskirts, and its mid-spans, with the unknown coefficients determined by minimizing the mean-square-error of the Bernoulli constant ($=0$).  For the highest wave with a corner crest of $120^{\circ}$, taking 15 unknown coefficients of the series expansion optimally selected has yielded results including $\alpha=a/h=0.8331990$ for the height, the Froude number $F=1.2908904$ for the wave speed, accurate to six decimals.  For almost highest waves, however, computations have been found to require noticeably more time and effort for results of comparable accuracy.  Further, for the so-called {\it dwarf long waves} defined for waves of amplitude $\alpha=a/h<10^{-2}$, including tsunamis commonly assessed as of amplitude $\alpha=O(10^{-3}-10^{-4})$ in the open ocean, this subject is still an open field in view of the findings[8] that the relative errors of the solution spread out much longer a stretch, lower the wave.  It can be strongly argued that in these challenging cases, knowing the singularity distribution of the solution outside the flow field should provide more accurate description of the wave properties, and hence improving the proficiency of the theoretical method.  It may further cast enlightening light on feasible mechanisms underlying wave instabilities, or bifurcation of the steady-state solution found.  Such pursuits can therefore be integrated with resolving the inverse problem as just described.  This is the principal objective of the present study.

To achieve this objective, it appears necessary first to extend the prevailing claims implied by Cauchy's integral formula, a fundamental theorem in theory of functions of a complex variable, $z=x+iy$, which we re-cite here as
\begin{subequations}\label{3 }
\be   J[f(z)]\equiv\frac{1}{2\pi i}\oint_C \frac{f(t)}{t-z}dt &=& f(z) \qquad~~~ (z\in {\cal D^+} -\mbox{open domain bounded by} ~C), \label{3a}\\
	&=& 0  	\qquad\qquad~ (z\in {\cal D^-} -\mbox{open domain outside} ~C),\label{3b} \nd
\end{subequations}
where the function $f(z)$, called Cauchy's function, is regular (holomorphic) (i.e. one-valued and continuously differentiable) inside the simply-connected open domain ${\cal D^+}$ bounded by contour $C$ and on $C$, $C$ being taken counter-clockwise, in positive sense.  The integral formula claims that the Cauchy integral, denoted by $J[f(z)]$, gives its value $f(z) ~\forall z\in {\cal D^+}$ and vanishes identically $~\forall z\in {\cal D^-}$ outside $C$.  Thus, $J[f(z)]$ has its value provided by the integral formula over the entire $z$-plane in terms of its values $f(t)$ on $C$, except $\forall z \in C$.  By this formula, the value of a continuous function $f(z)$ given on contour $C$ provides not only its local behavior, but also its global analytical properties over the domain $\cal D^+$ inside $C$, a distinguished property which is not shared by the function theory of a real variable.  Nevertheless, it has become evident that the crucial gap ($\forall z\in C$ being left uncovered by the formula) has kept the comprehensive properties of $f(z)$ and its integral formula from being thoroughly investigated and exposed.  Therefore, the foremost objective of the present study is to determine the value of the integral $ J[f(z)] ~\forall z\in C$ and to explore further associated developments.

For this purpose, we first introduce a new assumption that {\it $f(z)$ be $n$-times continously differentiable in ${\cal D}^+$ and in a neighborhood of contour $C$}.  Next, we conduct a limit procedure to let a point $z_+$ in the kernel $(t-z_+)^{-1}$ of $J[f(z)]\in {\cal D}^+$ and a point $z_-$ in the kernel $(t-z_-)^{-1} \in {\cal D}^-$ each tend from the ${\cal D}^\pm$ side to a common generic point $z_0 \in C$, respectively, by following two principles: 
\vskip 0.2mm 
\noindent (i) ~{\it Contour $C$ is deformed into $C^\pm$ with a small semicircle $C^{\pm}_\epsilon$ of radius $\epsilon$ centered at $z_0 \in C$ and indented onto the ${\cal D}^{\mp}$ side to let point $z_\pm (\in {\cal D}^\pm)$ reach $z_0$ without crossing $C^\pm$}, respectively.
\vskip 0.2mm 
\noindent (ii) ~{\it The value of Cauchy's integral remains intact in the limit as $z_{\pm}\rightarrow z_0$ and $\epsilon \rightarrow 0$}.\\
Fortunately, principles (i) and (ii) are just what are suitably needed in making two resulting relations (I) and (II) valid (see (\ref{9})).  Relation (I) proves the uniform continuity of $f(z)$ in the closed domain $\overline{{\cal D}^+} = [{\cal D}^+ +C]$, and (II) shows $f(z)$ at any $z$ on $C$ uniquely related to $f(t)~\forall t \in C$, as is fully delineated in \S3.  We ascribe these results of fundamental value to principle (i) in making indentation of $C$ possible and to principle (ii) in assuring the validity of the results, both under the new assumption in providing a neighborhood striding across $C$ for $f(z)$ to possess the analyticity, which is further extended to all the derivatives $J_n[f(z)]=d^n J[f(z)]/dz^n ~(n=1, 2, \cdots)$ of the integral formula.  The uniform continuity of $f(z)$ and  $f^{(n)}(z)=d^n f(z)/dz^n$ are readily applied in \S 4 to address the uniform convergence of their integral formulas.  From these fundamental findings there also follow the various integral properties of $f^{(n)}(z)$ and $J_n[f(z)]$ so produced.  

A new complement function $F(z)$ is introduced to be $C^n~\forall z\in {\cal D}^-$ outside $C$ and about $C$, and shown to share all the claims for $f(z)$ in complete analogy.  In \S 6, the uniform continuity found for Cauchy function $f(z)$ in domain $\overline{{\cal D}^+}$ and for the complement function $F(z)$ in $\overline{{\cal D}^-}= [{\cal D}^- +C]$ are jointly adapted and applied to special domains in the upper-, or lower-half $z$-planes and those inside or outside the unit circle $|z|=1$ to yield the generalized Hilbert transforms for these cases.

The central point is stressed that while the integral formula asserts $J[f(z)]\equiv 0 ~\forall z \in{\cal D}^-$, Cauchy's function $f(z)$ itself is nevertheless free to have any such singularity distributions in $\cal D^-$ as dictated only by its values $f(t)~\forall t \in C$.  The general behavior of Cauchy function $f(z)$ in $\cal D^-$ is illustrated in \S 7 as a {\it direct problem} (to which our {\it inverse problem} is its inverse) for various cases with  such singularities directly prescribed as poles, algebraic and logarithmic branches in $\cal D^-$.  

This study is extended in \S 8 and \S9 to include Plemelj's formulas for a line integral of a regular function along a Jordan arc without a double point ({\it not closed as a contour}) for applications to nonlinear wing theory.  It is concluded in \S 10 with expository discussions on the new results achieved here.   Finally, an unsolved inverse problem in regard to determining the singularity distribution of the Cauchy function in open domain $\cal D^-$ outside contour $C$ is presented for resolution as a conjecture.  This paper is also prepared with intent to enhance usage for self learning, research, and teaching.

\vskip 1mm
\noindent{\bf 2. ~The Cauchy function with Cauchy's integral theorem and integral formulas}.
\vskip 0.8mm
These fundamental theorems (cf. e.g.[1][6]) are recited in \S 2.1 to make this study self-contained.  
\vskip 1mm
\noindent{\bf 2.1. ~The Cauchy function $f(z)$ and the related classical theorems}.
\vskip 0.8mm
\noindent{\bf Definition 1. The Cauchy function.} ~{\it If analytic function $f(z)$ is regular in a simply-connected open domain ${\cal D}^+$ bounded by a simple contour $C$ and if $f'(z)\equiv df/dz$ is continuous on $C$, $f(z)$ is a Cauchy function.}  The open domain ${\cal D}^+$ plus contour $C$ plus the open domain ${\cal D}^-$ outside $C$ (including $z=\infty)$, i.e. $[{\cal D}^+ +C+{\cal D^-}]$ constitutes the entire $z$-plane. 
\vskip 0.5mm
\noindent{\bf Th.1. Cauchy's integral theorem.}~{\it If analytic function $f(z)$ is regular (holomorphic) inside and on a simple closed contour $C$, then the functional}
\be   I[f(z)]\equiv\oint_C f(z) dz =0 .  \label{1}\nd
We remark that Cauchy's theorem has been shown by Goursat to hold valid under a weaker assumption that $f'(z)$ exists (not necessarily being continuous) at all the points within and on $C$.  However, that the domain ${\cal D}^+$ be simply-connected is essential for the validity of (\ref{1}).  It leads directly to:

\vskip 1mm
\noindent{\bf ~Corollary C.1.} ~{\it If $f(z)$ is a regular function in a simply-connected domain ${\cal D}^+$, the integral
$$  \int_{z_0}^z f(t) dt = F(z)  \eqno(2a)$$
depends only on the end points $z_0$ and $z$, but not on the path between them in ${\cal D^+}$; and $F(z)$ is also a regular function in ${\cal D}^+$ such that $F'(z)=dF(z)/dz =f(z)$.}
\vskip 0.5mm
Conversely, if Corollary C.1 should be proved first, then Theorem Th.1 would follow as a corollary.  In fact, C.1 can be directly proved with $f(z)=u(x,y)+iv(x,y), f(z)dz =(udx-vdy)+i(vdx+udy)$, of which both the real and imaginary terms satisfy the condition of integrability in virtue of  
\be   \frac{\partial u}{\partial x} = \frac{\partial v}{\partial y}, \qquad  \frac{\partial u}{\partial y} = -\frac{\partial v}{\partial x} \qquad(\mbox{the Cauchy-Riemann equations}).  \label{2}\nd
For this proof for C.1, the same conditions as that invoked for Th.1 are necessary and sufficient.
\vskip 0.5mm
Next, we have Cauchy's integral formulas for $f(z)$ and its derivatives as follows.
\vskip 1mm
\noindent{\bf Theorem 2. Cauchy's integral formula.} This Theorem has already been re-cited in (1a,b).  There, both values of the integral $J[f(z)]$ for $z\in {\cal D}^{\pm}$ follow from (\ref{1}) of Theorem 1, since by Theorem 1, contour $C$ can be deformed to a small circle about any point $z$ in $C$ to yield (\ref{3a}) by the residue theorem, whereas for (\ref{3b}), the function $g(t,z)= f(t)/(t-z)$ is regular and satisfies (\ref{1}) $\forall~t\in C$ and $~\forall~z\in {\cal D}^-$.  
\vskip 0.8mm
\noindent{\bf Definition 2. The Cauchy integral and its derivatives.} ~{\it The integral defined in (1a,b) is called the {\it Cauchy integral}, here also called {\it the Cauchy functional} (to bear an integral operator connotation), to be denoted by $J[f(z)]$}, and its derivatives by $J_n[f(z)]\equiv d^nJ[f(z)]/dz^n$.
\vskip 0.5mm
\noindent{\bf Theorem 3. Derivatives of Cauchy's integral formula.} ~{\it Of Cauchy's integral formulas (\ref{3a})-(\ref{3b}), the functional $J[f(z)]$ has derivatives $J_n[f(z)]$ of all orders given by}
\begin{subequations}\label{4 a}
\be J_n[f(z)]\equiv\frac{n!}{2\pi i}\oint_C \frac{f(t)dt}{(t-z)^{n+1}} &=& f^{(n)}(z) \qquad ~(z\in {\cal D}^+;~n=1, 2, \cdots), \label{4a}\\
                 &=& 0  \qquad\qquad~~~ (z\in {\cal D}^-;~n=1, 2, \cdots).\label{4b} \nd \end{subequations}
These integral formulas are generally founded on the same basis as that for formulas (\ref{3a})-(\ref{3b}).  
\vskip 0.5mm
The foremost objective of this study is to determine, with proof, the values of $J_n[f(z)]~\forall z\in C ~(n=0,1,\cdots)$ so as to bring the closure of the integral formulas to completion.  For this general purpose, we first extend the original assumption invoked on the Cauchy function $f(z)$ as follows.
\vskip 0.5mm
\noindent{\bf 2.2. ~Generalization under new assumption}. ~{\it Here, we adopt a new assumption that $f(z)$ be $C^n$ ($n$-times continuously differentiable)~$\forall~z \in {\cal D}^+$ and in a neighborhood of contour $C$, $n$ being arbitrary}.  
\vskip 0.5mm
\noindent{\bf Definition 3. ~Generalized Cauchy function.}~The Cauchy function as just specified under the new assumption will be called the {\it generalized Cauchy function} when needed for its discreet distinction.
\vskip 0.5mm
\noindent{\bf Theorem 4. Derivatives of generalized Cauchy's integral formula.}~{\it With (\ref{4a})-(\ref{4b}) so extended, its integral can then be integrated by parts $m (\leq n)$ times, giving
\begin{subequations}\label{5 }
\be J_{(n,m)}[f(z)]\equiv \frac{(n-m)!}{2\pi i}\oint_C \frac{f^{(m)}(t)dt}{(t-z)^{n-m+1}} &=& f^{(n)}(z)\quad (z\in{\cal D}^+; ~m=0, 1, \cdots, n), \label{5a} \\
			&=& 0  \qquad\quad ~~(z\in {\cal D}^- ; ~n=1,2,\cdots), \label{5b} \nd \end{subequations}
in which every integrated term in each step, being single-valued on $C$, vanishes.} 

\vskip 1mm
\noindent{\bf ~Corollary C.4.} ~{\it In (\ref{5a}), $J_{(n,m)}[f(z)]$ are $(n+1)$ equivalent formulas for $f^{(n)}(z)$.  For $m=n$,}
\begin{subequations}\label{6}
\be  J_{(n,n)}[f(z)]\equiv\frac{1}{2\pi i}\oint_C \frac{f^{(n)}(t)}{t-z}dt &=& f^{(n)}(z) \qquad~~ (z\in {\cal D^+}; ~n=1,2,\cdots), \label{6a} \\
                &=& 0  \qquad\qquad\quad (z\in {\cal D}^- ;~n=1,2,\cdots).\label{6b} \nd \end{subequations}
This new result for the generalized Cauchy function $f(z)$ is very valuable, for writing $h(z)=f^{(n)}(z)$ renders (6a,b) identical in form with (1a,b), and can therefore be uniformly treated all together.
\vskip 0.5mm
\noindent{\bf Definition 4.~Complement function.} Interchanging the roles of ${\cal D}^{\pm}$ gives rise to the {\it complement function}, denoted by $F(z)$, defined in complete analogy in premise with Cauchy function $f(z)$, i.e. being $C^n~\forall ~z\in {\cal D}^{-}$ and in a neighborhood striding across contour $C$, including $z=\infty$ such that 
\be  F(z)=O(|z|^{-m}) \quad~~ (m\geq 2) \quad~~\mbox{as} \quad |z|\rightarrow \infty. \label{AM0} \nd
The integral of its integral formula, denoted by $J^-[F(z)]$, will be called the {\it complement functional}. 
\vskip 0.8mm
\noindent{\bf 2.3. ~The integral theorem and integral formula for complement function $F(z)$}.  ~For complement function $F(z)$, we first have
\vskip 0.5mm
\noindent{\bf Theorem 5. Integral Theorem for complement function $F(z)$}. {\it If complement function $F(z)$ is regular inside the simply-connected open domain ${\cal D}^-$ outside contour $C$, then the functional
\be  I^-[F(z)]\equiv \oint_{C^-} F(z)dz = 0,  \label{AM1} \nd
where contour $C^-$ is contour $C$ in opposite sense (each being in its own positive sense relative to ${\cal D}^{\pm}$). }
\vskip 0.5mm
\noindent{\it Proof:} ~Since $F(z)$ is regular in ${\cal D}^{-}$ and on $C$, the above contour integral around $C^-$ can be deformed, by Theorem 1, to that around $z=\infty$, without change in value, which vanishes by virtue of (\ref{AM0}).
\vskip 0.5mm 
To derive its integral formula, we set the origin $z=0$ inside contour $C$, apply the inverse conformal mapping $z=1/\zeta$ to flip the domains ${\cal D}^+$ and ${\cal D}^-$ in their images across the contour $C^*$ which is the image of contour $C$, and determine the integral formula for $F(z(\zeta))$ in the $\zeta$-plane by applying Cauchy's Theorem 2 to obtain, back in $z$ by inversion, the formula:
\begin{subequations}\label{AM2}
\be   J^-[F(z)]\equiv\frac{1}{2\pi i}\oint_{C^-} \frac{F(t)}{t-z}dt =\frac{-1}{2\pi i}\oint_{C} \frac{F(t)}{t-z}dt&=& F(z) \qquad~~~ (z\in {\cal D^-}), \label{AM2a}\\
	&=& 0  	\qquad\qquad~~ (z\in {\cal D^+}), \label{AM2b} \nd \end{subequations}
in which the negative sign for the integral around $C$ is due to the coincident contours $C^-$ and $C$ being opposite in sense.  More specifically, we apply the inverse conformal map $z=1/\zeta, t=1/\tau$ with respect to the origin $z=0$ which is set inside contour C, with $C$ mapped onto $C^*$ (also taken counter-clockwise, with ${\cal D^-}$ mapped onto ${\cal D^+_*}$ inside $C^*$ and ${\cal D^+}$ onto ${\cal D^-_*}$ outside $C^*$) so that 
\be  (2\pi i)J^-[F(z)]= \oint_{C^*} F\left(\frac{1}{\tau}\right)\frac{\zeta d\tau}{(\tau-\zeta)\tau} =\oint_{C^*} F\left(\frac{1}{\tau}\right)\left(\frac{1}{\tau-\zeta}-\frac{1}{\tau}\right)d\tau = \oint_{C^*} \frac{F(1/\tau)}{\tau-\zeta} d\tau \notag \nd
in which $\oint_{C^*} F(1/\tau)d\tau/\tau=2\pi iF(\infty)= 0$ since $F(\infty)$ (the residue at $\tau =0$) vanishes by (\ref{AM0}), whereas the final integral gives, by (1a,b) of Theorem 2, the resulting values in (\ref{AM2a})-(\ref{AM2b}).  We note that the integral for $J^-[F(z)]$ vanishes as $|z|\rightarrow \infty$ and (\ref{AM2a}) reduces to identity in this limit. 
\vskip 0.5mm
Extending (\ref{AM2a})-(\ref{AM2b}) to cover their derivatives in analogy with (1a,b) extended to (6a,b) then gives
\vskip 0.5mm
\noindent{\bf Theorem 6. Integral formula for complement function $F(z)$.} ~{\it The functional $J^-[F(z)]$ and its derivatives $J_n^-[F(z)]$ of a complement function $F(z)$ satisfy the integral formula} ($\forall n=0, 1, \cdots$):
\begin{subequations}\label{7}
\be   J^-_n[F(z)]\equiv\frac{1}{2\pi i}\oint_{C^-} \frac{F^{(n)}(t)}{t-z}dt=\frac{-1}{2\pi i}\oint_{C} \frac{F^{(n)}(t)}{t-z}dt &=& F^{(n)}(z) \qquad~~ (z\in {\cal D^-}), \label{7a}\\
	&=& 0  	\qquad\qquad~~~~ (z\in {\cal D^+}), \label{7b} \nd \end{subequations}
\vskip 0.5mm
Summing up this section, we note that that the values of $J_n[f(z)]$ and $J^-_n[F(z)]$ are now given in the entire $z$-plane for all $n = 0, 1, \cdots$ in entirety, except $\forall z \in C$.  We now pursue to determine their values $\forall z \in C$.  To proceed, we first resolve the apparent singularity of the kernel $(t-z)^{-1}$ for $\forall t\in C$ and also $\forall z\in C$ by using principles (i) and (ii) stated in \S1.
\vskip 0.5mm
\noindent{\bf 3. Uniform continuity of generalized Cauchy function $f(z)$ and complement function $F(z)$.} ~By Principle (i), we let point $z_\pm \in {\cal D}^{\pm}$ tend, respectively, to a generic point $z_0$ on $C$ by deforming $C$ into a closed contour $C^{\pm}=C_{\epsilon}^{\pm}+C_p$, where $C_{\epsilon}^{\pm}$ is a small semi-circle $|t -z_0|=\epsilon$ indented onto the ${\cal D}^{\mp}$-side of $C$ (so as to let $z_\pm$ reach $z_0 \in C$ without crossing $C^\pm$), and $C_p=C-C_{\epsilon}^{\pm}$ ~is preserved intact ~(see Fig. 1).  
\begin{figure}[hbt]
\begin{center}
\includegraphics[scale=0.64, trim=16mm 16mm 12mm 12mm]{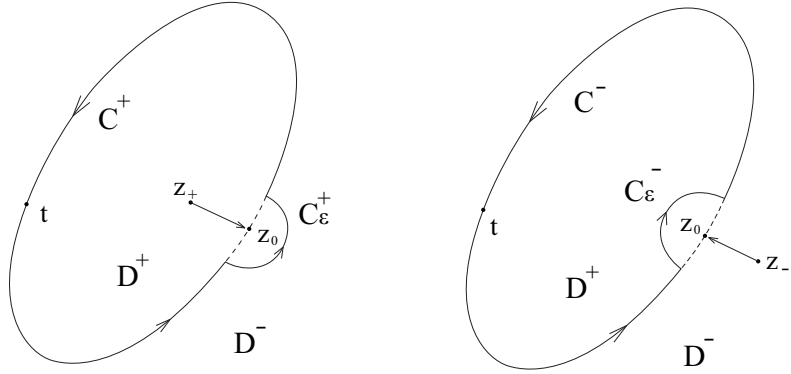}
\caption{\small  A Cauchy integral $\oint_C f(t)dt/(t-z)$ around a simple closed contour $C$ enclosing domain ${\cal D}^+$ and excluding domain ${\cal D}^{-}$ has its contour deformed into $C^{\pm}$ with only a small semicircle $C_{\epsilon}^{\pm}$ indented onto the ${\cal D}^{\mp}$ side, centered at $z_0\in C$ so as to let point $z_\pm\in {\cal D}^{\pm}$ tend, respectively, to $z_0\in C$ without crossing $C^{\pm}$.}
\end{center}
\end{figure}
\vskip -2mm 
\noindent Since as $z_\pm$ reaches $z_0$, $z_0$ on $C$ is never crossed by either $C^+$ or $C^-$, Principle (ii) is observed.  Hence 
\be & &I^{\pm}(z_0)\equiv \oint_{C^{\pm}}\frac{f(t)}{t -z_0}dt=\left\{\int_{C_{\epsilon}^{\pm}}+ \int_{C_p}\right\}\frac{f(t)}{t -z_0}dt = I_{\epsilon}^{\pm}(z_0)+I_{p}(z_0) \qquad\qquad (\epsilon> 0), \notag\\
 & &I_{\epsilon}^{\pm}(z_0)= f(z_0)\int_{C_{\epsilon}^{\pm}}\frac{dt}{t -z_0}+ J(z_0)\rightarrow \pm\pi i f(z_0), \quad J(z_0)=\int_{C_{\epsilon}^{\pm}}\frac{f(t)-f(z_0)}{t -z_0}dt~\rightarrow 0 \quad~  (\epsilon\rightarrow 0), \notag\nd
where the limit for $I_{\epsilon}^{\pm}(z_0)$ comes with $t$ moving in the $\pm$ sense on $C_{\epsilon}^{\pm}$, whilst that for $J(z_0)$ results from the argument that for given $\epsilon> 0$, $\exists ~\delta (\epsilon,z_0) \ni |(f(t)-f(z_0)) /(t -z_0)-f'(z_0)|<\delta, \forall~|t -z_0|<\epsilon$, hence $|J(z_0)|<(|f'(z_0)|+\delta)\ell_\epsilon ~(C_{\epsilon}^{\pm}$ being of length $\ell_\epsilon=\pi \epsilon)\rightarrow 0$ as $\epsilon \rightarrow 0$.  In addition, as $\epsilon \rightarrow 0$, the integral over $C_p$ assumes Cauchy's principal value defined by
\be  I_{p}(z_0)= \int_{C_{p}}\frac{f(t)}{t -z_0}dt \rightarrow \lim_{\epsilon\rightarrow 0}\left\{\int_{t_0}^{z_0-\epsilon}+\int_{z_0+\epsilon}^{t_0}\right\} \frac{f(t)dt}{t-z_0} = {\cal P}\oint_C \frac{f(t)dt}{t-z_0} \qquad (z_0 ~\mbox{on}~ C),\notag \nd
where ${t_0} (\neq z_0)$ is any point on $C$, $z_0\pm\epsilon$ both lying on the regular curve $C$, the symbol ${\cal P}$ signifying Cauchy's principal value being often omitted by convention as understood.  

Finally, as $z_+\rightarrow z_0$ and $\epsilon\rightarrow 0$, (\ref{3a}) yields $J[f(z_0)]= (I_{\epsilon}^+ +I_p)/(2\pi i)\rightarrow f^+(z_0)$, a value of the limit of $f(z)$ as $z_+\rightarrow z_0$, yet to be determined.  And as $z_-\rightarrow z_0, \epsilon\rightarrow 0$, (\ref{3b}) gives $J[f(z_0)]=(I_{\epsilon}^- +I_p)/(2\pi i)\rightarrow 0$, since $J[f(z_0)]$ of (\ref{3b}) remains zero by Principle (ii).  Thus, (\ref{3a})-(\ref{3b}) become
\begin{subequations}\label{8}
\be  f^+(z_0)&=& \lim_{z_+\rightarrow z_0}\frac{1}{2\pi i}\oint_C \frac{f(t)}{t-z_+}dt=~\frac{1}{2}f(z_0)+\frac{1}{2\pi i}{\cal P}\oint_C \frac{f(t)}{t-z_0}dt ~\quad ~~(z_0 ~\mbox{on}~ C),\label{8a}\\
           0 &=& \lim_{z_-\rightarrow z_0}\frac{1}{2\pi i}\oint_C \frac{f(t)}{t-z_-}dt= -\frac{1}{2}f(z_0)+\frac{1}{2\pi i}{\cal P}\oint_C \frac{f(t)}{t-z_0}dt\quad ~~(z_0 ~\mbox{on}~ C),\label{8b}\nd \end{subequations}
of which the sum and difference determine two relations of vital importance as
\be  (I):~~~f^+(z)=f(z); \quad~~ (II):~~~f(z)=\frac{1}{\pi i}{\cal P}\oint_{C}\frac{f(t)}{t-z}dt \qquad (z ~\mbox{on}~ C), \label{9} \nd
in which the suffix of $z_0$ is omitted for all $z$ on $C$.  Here, the first relation (I), $f^+(z)=f(z)$, shows that the limit $f^+(z)$ of $f(z)$ reached from the ${\cal D}^+$ side is equal to the original $f(z)$ prescribed $ \forall~z\in C$, therefore proves the {\it uniform continuity of $f(z)$ in the closed domain} $\overline{{\cal D}^+}=[{\cal D}^+ +C]$.  Relation (II) gives the value to $f(z)$ for any $~z\in C$ in terms of its values $f(t)$ over $C$.  The integral representing $f(z)$ for $z\in C$ in (II) is noted to differ from that of $f(z)$ for $z\notin C$ in (\ref{3a}) by a factor of $2$, with the new integral assuming its principal value.  This completes the proof of the above two key relations by finding the limits (\ref{8a})-(\ref{8b}) of $J[f(z)]$ of {\it both} (\ref{3a}) and (\ref{3b}).  By (\ref{9}), we therefore have proved
\vskip 1mm
\noindent{\bf Theorem 7. Uniform continuity of Cauchy function $f(z)$ in closed domain $\overline{{\cal D}^+}$.} {\it If $f(z)$ satisfies Cauchy's integral formula (\ref{3a}) in open domain $\cal D^+$ bounded by contour $C$ and also (\ref{3b}) in open domain $\cal D^-$ outside $C$, then $f(z)$ is uniformly continuous within the closed domain $\overline{{\cal D}^+}=[{\cal D}^+ +C]$.}
\vskip 0.5mm
Since the extended formulas (6a,b) are identical in form with (1a,b), we have also proved 
\vskip 0.5mm
\noindent{\bf Theorem 8. Uniform continuity of the derivatives $f^{(n)}(z)$ in closed domain $\overline{{\cal D}^+}$.} {\it If $f(z)$ satisfies (\ref{3a}) in open domain $\cal D^+$ and also (\ref{3b}) in open domain $\cal D^-$, and if $f^{(n)}(z)$ is continuous on $C$, then $f^{(n)}(z)$ is uniformly continuous within the closed domain $\overline{{\cal D}^+}=[{\cal D}^+ +C]$,}
\be  f^{(n)+}(z)=f^{(n)}(z); \quad~~ f^{(n)}(z)=\frac{1}{\pi i}{\cal P}\oint_{C}\frac{f^{(n)}(t)}{t-z}dt  \qquad ~~(z ~\mbox{on}~ C;~n=1, 2, \cdots).  \label{10} \nd
\vskip 1mm
\noindent{\bf 4. Uniform convergence of the integral formulas and resulting integral properties.}
\vskip 0.8mm
From Theorems 7 and 8 there readily follow other important consequences:
\vskip 0.8mm
\noindent{\bf Theorem 9. Uniform convergence of the integral formulas.} ~{\it  The integral formulas of $J_n[f(z)] ~\forall z\in {\cal D}^+$ (or $~\forall z\in {\cal D}^-$) converge uniformly in closed domain $\overline{{\cal D}^+}$ (or $\overline{{\cal D}^-}$), respectively,  for $n=0,1,2 \cdots$.} 
\vskip 0.8mm
\noindent{\it Proof}: Rewrite (1a,b) and (6a,b) for $n=0,1,2 \cdots$ as
\be  g_n(z)\equiv J_n[f(z)]- f^{(n)}(z)=0 ~~~(z \in {\cal D}^+); \qquad G_n(z)\equiv J_n[f(z)] =0 ~~~(z \in {\cal D}^-), \label{AM5} \nd
where $J_n[f(z)]$ stands for the integral it represents for brevity.  In (\ref{AM5}), let $z_{\pm} \in {\cal D}^\pm$ tend, respectively, to a generic point $z_0 \in C$, with contour $C$ indented into $C_\pm$ as delineated before (see Fig. 1), thus yielding (\ref{8a})-(\ref{8b}) for the values of the integral $J_n[f(z)]$ in the limit.  These two limiting equations, upon using the two relations in (\ref{9}), then become
\be  g_n(z_0)=G_n(z_0)= -\frac{1}{2}f^{(n)}(z_0)+\frac{1}{2\pi i}{\cal P}\oint_{C}\frac{f^{(n)}(t)}{t-z_0}dt =0 \quad ~~(z_0 ~\in~ C;~n=0,1, 2, \cdots).  \label{AM6} \nd
Therefore, (\ref{AM5}) and (\ref{AM6}) jointly state that for $z_0 ~\in~ C, ~n=0,1, 2, \cdots$,
\be g_n(z)= J_n[f(z)]- f^{(n)}(z)= g_n(z_0)=0 ~~~(z \in \overline{{\cal D}^+});~~~~ G_n(z)\equiv J_n[f(z)] =G_n(z_0)=0 ~~~(z \in \overline{{\cal D}^-}), \label{AM7} \nd 
which clearly proves the uniform convergence of $g(z)$ in closed domain $\overline{{\cal D}^+}$ and that of $G(z)$ in closed domain $\overline{{\cal D}^-}=[{\cal D}^- + C^-]$, as was to be proved.
\vskip 0.5mm
Concerning the uniform continuity of complement function $F(z)$ in analogy with $f(z)$, we have noted that the integral $J^-[F(z(\zeta))]$ in the $\zeta$-plane given by the $z=1/\zeta$ map is identical in form with the integral of $J[f(z)]$ in the $z$-plane of (\ref{3a}). Hence we can directly apply the limiting formula (\ref{8a})-(\ref{8b}) (for the integral formula of $J[f(z)]$) to the integral formula of $J^-[F(z(\zeta))]$ in the $\zeta$-plane, giving, after inverting back to the $z$-plane, the result in analogy with (\ref{10}) as
\be F^{(n)-}(z)=F^{(n)}(z); \quad~~F^{(n)}(z)=\frac{1}{\pi i}{\cal P}\oint_{C^-}\frac{F^{(n)}(t)}{t-z}dt =\frac{-1}{\pi i}{\cal P}\oint_{C}\frac{F^{(n)}(t)}{t-z}dt  \quad ~~(z ~\mbox{on}~ C),  \label{11} \nd 
for $n=0, 1, \cdots$, where $F^{(n)-}(z)$ is the limit of $F^{(n)}(z)$ as $z$ reaches $C$ from the ${\cal D}^-$ side. Whence
\vskip 0.8mm
\noindent{\bf Theorem 10. Uniform continuity of complement function $F(z)$ and its derivatives $F^{(n)}(z)$ in closed domain $\overline{{\cal D}^-}$ and uniform convergence of their integral formulas.} {\it If $F(z)$ satisfies (\ref{7a})   $\forall z \in {\cal D}^-$, and (\ref{7b}) $\forall z \in {\cal D}^+$, then $F(z)$ and its derivatives $F^{(n)}(z)$ are uniformly continuous in closed domain $\overline{{\cal D}^-}$, and their integral formulas are uniformly convergent in analogy to that for $f(z)$}. 
\vskip 0.8mm
We note that in (\ref{9})-(\ref{10}) and (\ref{11}) there exist other distinct functionals of $f(z)$ and $F(z)$ as
\be  K_n[f(z)]\equiv \frac{1}{\pi i}{\cal P}\oint_{C}\frac{f^{(n)}(t)}{t-z}dt = f^{(n)}(z), ~\quad K_n^-[F(z)]\equiv \frac{-1}{\pi i}{\cal P}\oint_{C}\frac{F^{(n)}(t)}{t-z}dt = F^{(n)}(z) \quad (z ~\mbox{on}~ C), \label{12} \nd
for $n= 0, 1, 2, \cdots$, with $z$ strictly lying on contour $C$.  These functionals have merits of their own.
\vskip 0.8mm
Regarding uniform continuity and uniform convergence, Theorems 7 to 9 for $f(z)$ and Theorem 10 for $F(z)$ are of fundamental importance because to them all the general theorems on uniform continuity and uniform convergence (see e.g. [1]) then hold.  The six limiting equations in (\ref{9})-(\ref{10}) and (\ref{11}) thus set the foundation established here for further developments, some to follow next.
\vskip 0.8mm
\noindent{\bf Th.11. Integral theorem of functionals $J_n[f(z)]$.} {\it The contour integrals of functionals $J_n[f(z)]$ defined by limiting equations (\ref{9})-(\ref{10}) and that of $J_n^-[F(z)]$ by (\ref{11}) along contour $C$ all vanish},
\begin{subequations}\label{13}
\be \oint_C J_n[f(z)]dz =\oint_C f^{(n)}(z)dz =0 && (n=0, 1, 2, \cdots), \label{13a}\\
\oint_C J_n^-[F(z)]dz =\oint_C F^{(n)}(z)dz =0 && (n=0, 1, 2, \cdots) . \label{13b}\nd \end{subequations}
{\it Proof:} ~In fact, (\ref{13a}) follows from the uniform convergence relation (\ref{AM7}) and the uniform continuity of $f^{(n)}(z)$ within closed domain $\overline{{\cal D}^+}$, and likewise (\ref{13b}) holds within closed domain $\overline{{\cal D}^-}$ by analogy.  
\vskip 0.8mm
\noindent{\bf Th.12. Integral theorem of functionals $K_n[f(z)]$ and $K_n^-[F(z)]$.} {\it The contour integrals of functionals $K_n[f(z)]$ and $K_n^-[F(z)]$ of (\ref{12}) along contour $C$ all vanish,}
\begin{subequations}\label{14}
\be \oint_C K_n[f(z)]dz =\oint_C f^{(n)}(z)dz =0 && (n=0, 1, 2, \cdots), \label{14a}\\
\oint_C K_n^-[F(z)]dz =\oint_C F^{(n)}(z)dz =0 && (n=0, 1, 2, \cdots) . \label{14b}\nd \end{subequations}
\vskip 0.8mm
\noindent{\it Proof}: ~Rather than seeking a proof by analogy with that for (\ref{13}), we give a direct proof with $z$ strictly lying on contour $C$.  Integrating (\ref{12}) over $C$ and interchanging the order of integration yields
\be  \oint_C f^{(n)}(z)dz = \frac{1}{\pi i}{\cal P}\oint_C dz \oint_C \frac{f^{(n)}(t)}{t-z}dt =\frac{1}{\pi i}{\cal P}\oint_C f^{(n)}(t)dt\oint_C \frac{dz}{t-z}= -\oint_C f^{(n)}(t)dt , \nonumber\nd
since ${\cal P}\oint_C dz/(t-z) =-\pi i ~\forall t\in C$.  To prove this, we consider its integral over a closed contour $C^+=C_{\epsilon}^+ +C_p$ as shown in Fig. 1, with the $C_{\epsilon}^+$ semi-circle centered at point $t$ on $C$, so that
\be  \oint_{C^+} \frac{dz}{z-t} = 2\pi i = \lim_{\epsilon\rightarrow 0}\left\{\int_{C_{\epsilon}^{+}}+ \int_{C_p}\right\}\frac{dz}{z-t} = \pi i + {\cal P}\oint_{C}\frac{dz}{z-t} \quad\rightarrow\quad {\cal P}\oint_C \frac{dz}{t-z}= -\pi i, \nonumber \nd
hence (\ref{14a}) results, for $\oint_C f(z)dz$ is found to be negative of itself.  Likewise, (\ref{14b}) holds by analogy.    
\vskip 0.8mm
We may note that the integral formula (\ref{4a}) for $f^{(n)}(z)$ can be attained by direct differentiation of formula (\ref{3a}) for $f(z)$ {\it under the integral sign}.  This operation can be justified as follows.
\vskip 0.8mm
\noindent{\bf Theorem 13. Differentiation of functional $J[f(z)]$.} {\it The derivatives $J_n[f(z)] ~(n=1, 2,\cdots)$ of functional $J[f(z)]$ can be derived by direct differentiation of $J[f(z)]$ under the integral sign.} 
\vskip 1mm  This is valid due to $J[f(z)]$ (before differentiation) and $J_n[f(z)]$ (after differentiation) being all uniformly continuous in any closed sub-domain $\overline{{\cal D}_s}$ inside the open domain $\cal D^+$ and also in the closed domain $\overline{{\cal D}^+}$ in virtue of $J_n[f(z)]=J_{(n,n)}[f(z)]$ of (\ref{6a}) being identical in form with (\ref{3a}).  

\vskip 1mm
\noindent{\bf 5. Mean value of regular functions.} ~Formulas (\ref{3a}) and (\ref{6a}) provide a mean value of $f(z)$ and $f^{(n)}(z)$ for a circular contour in particular.  The discussions here will be addressed only for $f(z)$ as the corresponding results for complement function $F(z)$ can be implied by analogy.
\vskip 1mm
\noindent{\bf Th.14. Mean-value theorem.} {\it If each of $f^{(n)}(z)$ exists and is continuous inside and on a circular contour $C$, $t -z =r \exp (i\theta)$ ~(of arbitrary radius $r$), formulas (\ref{3a}) and (\ref{6a}) become}
\be   f^{(n)}(z)=\frac{1}{2\pi}\int_0^{2\pi} f^{(n)}(z+r e^{i\theta}) d\theta \qquad (n=0, 1, 2, \cdots), \label{15}\nd
asserting that the value of a regular function $f^{(n)}(z)$ at the center of a circle $C$ is equal to the mean of its values on $C$.  As $r\rightarrow 0$, (\ref{15}), being homogeneous in $f^{(n)}$, becomes an identity.  

In general, let $|f^{(m)}(t)|\leq M_m(R,z)$ be the upper bound on $|t -z|=R$, then by (\ref{5a}),
\be  |f^{(n)}(z)|\leq (n-m)!~R^{-(n-m)}M_m(R,z)\qquad (R=|t -z|, ~m=0,1,\cdots n; ~n=0,1,2,\cdots). \label{16} \nd
To explore dependence of $M_m(R,z)$ on $R=|t -z|$ and $z$, let us first consider a class of function $f(z)$ that is regular in an open domain ${\cal D}^+$ and hence possesses a Taylor series $f(z)= \Sigma_{n=0}~c_nz^n$, convergent absolutely inside a circle $|z|=R$ lying within ${\cal D}^+$.  Let $|f(z)|\leq M=\mbox{max}_{|z|=R}|f(z)|$ for $|z|\leq R$, then
\be   |c_n|R^n\leq M (R) \qquad (\mbox{Cauchy's inequality}),  \label{17} \nd
which in turn provides for $|f^{(n)}(0)|$, in virtue of $f^{(n)}(0)=n!~c_n$, the upper bounds 
\be  |f^{(n)}(0)|\leq n!R^{-n}M (R)\qquad (n=0, 1, 2, \cdots).  \label{18} \nd
As a proof for (\ref{17}) and (\ref{18}), we deduce from (\ref{4a}), for $|t|=R$, the relation
\be  |f^{(n)}(0)|\leq \frac{n!}{2\pi}\oint_C \frac{M}{R^{n+1}}|dt| =n!R^{-n}M (R),  \nonumber \nd
in agreement with (\ref{16}) for $z=0$ and $m=0$.  An immediate consequence of (\ref{18}) is the following.
\vskip 0.8mm
\noindent{\bf Th.15. Liouville's theorem.} {\it If $f(z)$ is analytic and bounded for all finite $z$, it is a constant.}
\vskip 1mm In fact, if $|f(z)|\leq M$ for all $z$, $|f^{(n)}(0)|\rightarrow 0$ as $R=|z|\rightarrow \infty$ by (\ref{18}) for all $n\geq 1$, leaving only the $n=0$ term to give $f(z)=c_0$, a constant.  By an extension in scope, we have
\vskip 0.8mm
\noindent{\bf Th.16. Polynomial theorem.} {\it If $f^{(m)}(z)$ is analytic and $|f^{(m)}(z)|\leq M_m ~(\mbox{const}.>0) ~\forall z$, $f(z)$ is a polynomial of degree $m$.}
\vskip 0.8mm First, by Liouville's theorem, $f^{(m)}(z)$ is a constant, hence by (\ref{16}), $|f^{(n)}(0)|\rightarrow 0$ as $R\rightarrow \infty ~\forall n\geq m+1$, as was to prove.  This theorem may be called the {\it extended Liouville's theorem.}

\vskip 1mm
\noindent{\bf 6. The generalized Hilbert transforms.} ~We next pursue whether there exists an integral analog of the Cauchy-Riemann differential relations (\ref{2}) between the conjugate functions $u$ and $v$ of an analytic function $f(z)=u(x,y)+iv(x,y)$.  This leads to Hilbert's integral transform we now discuss.

\noindent{\bf 6.1. The Hilbert transform.} ~We consider first a class of analytic function $f(z)$ which is regular in the upper half $z$-plane for $Im~z\geq 0$, and vanishes as $|z|\rightarrow \infty$ uniformly in $0\leq \arg z\leq \pi$; then by formula (\ref{3a})-(\ref{3b}) for this $f(z)$ we take $C$ along the upper semicircular contour $C_u =C_x (-R\leq x\leq R)+C_R^+ (z=R e^{i\theta}, R=$ const., $0\leq\theta\leq\pi)$, the integral on $C_R^+\rightarrow 0$ as $R\rightarrow \infty$, giving
\be  \oint_{C_u}\frac{f(z)}{z-\zeta}~dz =\int_{-\infty}^{\infty}\frac{f(x)}{x-\zeta}~dx &=& 2\pi i f(\zeta) \qquad ~~~(Im~\zeta >0), \nonumber \\
		&=& 0 \qquad\qquad \qquad (Im~\zeta <0). \nonumber\nd
The limit of this equation as $\zeta\rightarrow \xi$ (a point on the real $\zeta$-axis), from above or from below, has been obtained for arbitrary contour in (\ref{9}) which can be adapted to the present geometry to give
\be f^+(\xi)=f(\xi); \quad~~ f(\xi)=\frac{1}{\pi i}{\cal P}\int_{-\infty}^\infty\frac{f(x)dx}{x-\xi}\qquad~(-\infty<\xi<\infty).\label{21}\nd
This shows that by Theorem 7, $f(z)$ is uniformly continuous in the closed domain $\overline{{\cal D}_u^+}:(0\leq |z|\leq R<\infty, 0\leq \arg z\leq \pi)$.  Hence substituting $f(x)=u(x)+iv(x), ~f(\xi)=u(\xi)+iv(\xi)$ in (\ref{21}), with $(u, v)\in {\cal C}^1 (-\infty<x<\infty)$ being understood, yields for the real and imaginary parts as
\be  u(\xi)=H[v(x)]=\frac{1}{\pi}{\cal P}\int_{-\infty}^\infty\frac{v(x)dx}{x-\xi}, \qquad
 v(x)=H^{-1}[u(\xi)]= \frac{-1}{\pi}{\cal P}\int_{-\infty}^\infty\frac{u(\xi)d\xi}{\xi-x} . \label{22} \nd
This pair of reciprocal integral relations, known as the {\it Hilbert transform}, is due to David Hilbert (1862-1943), with $H$ denoting the transform and $H^{-1}$ the inverse transform.  In relations (\ref{22}), $u(x)$ is said to be conjugate to $v(x)$; the relationship is {\it skew-reciprocal}, i.e. reciprocal apart from a minus sign, e.g. $-v(x)$ is conjugate to $u(x)$.  For a function regular for $Im~z\geq 0$, e.g. $f(z)=e^{iz},~ e^{ix}=\cos x+i\sin x$, we have $\cos\xi=H[\sin x]$, whereas by inversion, $\sin (x)=H^{-1}[\cos \xi]~=-H[\cos \xi]$, in skew reciprocity.  Now, by direct symbolic substitutions of the two operator equations (\ref{22}), we obtain the relations 
\be   H^{-1}H[v(x)]=v(x), \qquad HH^{-1}[u(\xi)]=u(\xi), \quad\longrightarrow\quad H^{-1}H = HH^{-1}= 1 \quad(\mbox{unity operator}). \label{23} \nd
This can be shown for specific $u(x)$ or $v(x)$ by consecutive evaluation of the integrals as exemplified here, whereas  showing this for arbitrary $u(x)$ or $v(x)$ will involve interchanging the order of integrations involving product of two Cauchy kernels.  For dealing with such a general case, it is essential to have
\vskip 0.8mm
\noindent{\bf The Poincar$\acute{e}$-Bertrand formula:}
\be  \int_L\frac{dt'}{t'-x}\int_L\frac{f(t,t')}{t-t'}dt = \int_L dt\int_L\frac{f(t,t')dt'}{(t'-x)(t-t')} - \pi^2 f(x,x) \qquad  (x \in L), \label{36} \nd
where $L$ is a regular Jordan arc, assumed finite (or infinite) in length, with end-points at $t=a$ and $t=b$ and without double point, the integration variable $t$ moves from $a$ to $b$, and function $f(t,t')$ is assumed regular in a neighborhood of the entire line $L$, while each of the integrals assumes its own principal value, here with the symbol $\cal P$ omitted as understood by convention.  For its proof we refer to the literature (e.g. Muskhelishvili[4]) and a hint given for the general formula (\ref{36}).  To illustrate application of the Poincar$\acute{e}$-Bertrand formula (\ref{36}) to Cauchy integrals, we return to (\ref{23}) for a proof in general as follows.
\vskip 0.8mm
\noindent{\it Example 1.} ~Consider the formula $H^{-1}H[v(x)]= v(x)$ for arbitrary $v(x)$ being operated by the Hilbert transform and its inversion in succession so that
\be  H^{-1}H[v(x)]= \frac{-1}{\pi^2}\int_{-\infty}^\infty \frac{dt'}{t'-x}\int_{-\infty}^\infty \frac{v(t)dt}{t-t'} =  v(x)+\frac{1}{\pi^2}\int_{-\infty}^\infty \frac{v(t)dt}{x-t}\int_{-\infty}^\infty (\frac{1}{t'-x}-\frac{1}{t'-t})dt' = v(x),  \notag\nd
where the second equality results from applying the Poincar$\acute{e}$-Bertrand formula and the last integral vanishes since ${\cal P}\int dt'/(t'-x) =0$.  Similarly, we can show that $HH^{-1}[u(x)]=u(x)$, thus providing a proof of the relations in (\ref{23}) previously implied by substitutions of the two operator equations.

Indeed, this also shows that when the two key relations in (\ref{9}) can be applied to a function $f(z)=u(x,y)+iv(x,y)$, regular in a certain domain, to result in a skew-reciprocal pair of transform equations between its conjugate functions $u$ and $v$ (like that in the present case for the upper-half $z$-plane and three more to follow), the transform relations arrive automatically, with no need to prove that $u$ and $v$ are conjugate functions (like in some other approaches[7]), for the proof is already imbedded in (\ref{9}).  Also owing to (\ref{23}), we can assert that if $u(\xi)=H[v(x)]$ is regarded as a singular integral equation for $v(x)$ with $u(\xi)$ given (being Hilbert transformable), its solution is $v(x)=H^{-1}[u(\xi)]$, and vice versa.
\vskip 0.8mm
In this case, if we find the Hilbert transform $u(\xi)=H[v(x)]$ of $v(x), C^1~\forall x(-\infty,\infty)$ to form a complex function $f(x)=u(x)+iv(x)$, and have it analytically continued into $f(z)$ over the entire $z$-plane, then, by implication of the analysis underlying (\ref{21}), $f(z)$ must be analytic and regular in the upper half $z$-plane, and further, by Theorem 7, be uniformly continuous in the closed domain $\overline{{\cal D}_u^+}$.  
\vskip 0.8mm
\noindent{\it Example 2.} ~If we take $v(x)=\cos x$, then by (\ref{22}), $u(\xi)=H[v(x)]=H[\cos x]=-\sin \xi$, hence $f(x)=u(x)+iv(x)=i(\cos x +i\sin x)=i e^{ix}$, giving the analytically continued function $f(z)=ie^{iz}$, which is regular in the upper half $z$-plane, but singular in the complementary lower-half $z$-plane.

\vskip 1mm
\noindent{\bf 6.2. The complementary Hilbert transform.} ~On the contrary, if $F(z)$ is regular in the lower half $z$-plane, then we have $F(z)$, regular in ${\cal D}^- \{z|Im ~z \leq 0\}$, satisfying, by (\ref{11}), the relations 
\be F^-(\xi)=F(\xi); \quad~~ F(\xi)= \frac{-1}{\pi i}{\cal P}\int_{-\infty}^\infty\frac{F(x)dx}{x-\xi}\qquad~(-\infty<\xi<\infty),\label{24}\nd 
where $F^-(\xi)$ is the limit of $F(z)$ as a point $z ~(Im ~z <0)$ tends from below to reach a point $\xi$ on the real $z$-axis, so that (\ref{24}) differs from the corresponding relation (\ref{21}) only by a minus sign of the integral by virtue of (\ref{11}).  Thus, substituting $F(x)=U(x)+iV(x),~F(\xi)=U(\xi)+iV(\xi)$ in (\ref{24}) yields
\begin{subequations}\label{25 }
\be  U(\xi)=\overline{H}[V(x)]=\frac{-1}{\pi}{\cal P}\int_{-\infty}^\infty\frac{V(x)dx}{x-\xi}, & &
 V(x)=\overline{H}^{-1}[U(\xi)]= \frac{1}{\pi}{\cal P}\int_{-\infty}^\infty\frac{U(\xi)d\xi}{\xi-x} , \label{25a} \\
   \longrightarrow \quad \overline{H}[V(x)]= -H[V(x)]= H^{-1}[V(x)], & & \overline{H}^{~-1}[U(\xi)]=H[U(\xi)]=-H^{-1}[U(\xi)].  \label{25b} \nd \end{subequations}
The pair of reciprocal integral relations in (\ref{25a}), designated by $\overline{H}[\cdot]$ and its inverse by $\overline{H}^{-1}[\cdot]$, may be called the {\it complementary Hilbert transform}; it is related to the Hilbert transform by (\ref{25b}).  
\vskip 1mm
\noindent{\it Example 3.} ~Given $v(x)=-(x^2+1)^{-1}$, its Hilbert transform is given, after some algebra, by
\be   u(\xi)=\frac{-1}{\pi}{\cal P}\int_{-\infty}^{\infty}\frac{dx}{(x^2+1)(x-\xi)}= \frac{\xi}{\xi^2+1}~, \quad f(x)=u(x)+iv(x)=\frac{1}{x+i} ~~\rightarrow~~ f(z)=\frac{1}{z+i} \nonumber \nd
which is regular in the upper half $z$-plane, but has a simple pole at $z=-i$ in the lower half $z$-plane.  On the other hand, the complementary transform of $V(x)=-(x^2+1)^{-1}$ gives, by (\ref{25a}),
\be   U(\xi)= \overline{H}[V(x)]= -H[V(x)]= \frac{-\xi}{\xi^2+1} \quad\rightarrow\quad F(x)= U+iV=\frac{-1}{x-i} ~~\rightarrow~~F(z)= \frac{-1}{z-i} , \nonumber \nd
which is regular in the lower half $z$-plane, but has a simple pole at $z=i$ in the upper half $z$-plane.  Thus a given $v(x)$ has been demonstrated to generate an analytic function $f(z)$ (or $F(z)$) which is regular in the upper (or lower) half $z$-plane by applying the Hilbert (or the complementary Hilbert) transform.  
\vskip 0.8mm
\noindent{\bf 6.3. The circular Hilbert transform.} We next consider function $f(z)$ which is regular in open domain ${\cal D}_c^+$ of a unit disc: $|z|\leq 1$ and $f(z) \in C^1$ on $|z|=1$, which we take for the contour $C$ in the general formula (\ref{9}) with  both points $z=e^{i\theta}$ and $t=e^{\i\phi}$ on $C ~(|t|=1$), giving, for ($-\pi\leq\theta\leq \pi$),
\begin{subequations}\label{26 }
\be && \qquad\qquad f(e^{i\theta}) = \frac{1}{\pi}{\cal P}\int_{-\pi}^\pi f(e^{\i\phi}) \frac{e^{i\phi} d\phi}{e^{i\phi}-e^{i\theta}} = \frac{1}{2\pi}{\cal P}\int_{-\pi}^\pi (1-i\cot \frac{\phi -\theta}{2})f(e^{\i\phi})d\phi, \label{26a} \\
&& u(\theta)= \frac{\cal P}{2\pi}\int_{-\pi}^\pi \left(u(\phi)+v(\phi)\cot \frac{\phi -\theta}{2}\right)d\phi, \quad~
     v(\phi) = \frac{\cal P}{2\pi}\int_{-\pi}^\pi \left(v(\theta)-u(\theta)\cot \frac{\theta-\phi}{2}\right)d\theta, \qquad~~\label{26b} \nd \end{subequations}
which results from separating the real and imaginary parts in (\ref{26a}) with $f(e^{i\theta})=u(\theta) +iv(\theta)$, 
with $u(-\pi)=u(\pi)$ and $v(-\pi)=v(\pi)$ understood.  Concerning the mixed functions in the integrals, we notice that in this case, (\ref{13}) of Th.11 or (\ref{14}) of Th.12 reduces to an integral of a single variable as
\be  \int_{-\pi}^\pi f(e^{i\theta}) d\theta =\int_{-\pi}^\pi \{u(\theta)+i v(\theta)\} d\theta =0. \label{27} \nd
This intrinsic normalization condition can be adopted to resolve discrepancies between some similar yet differing published expressions for conjugate equations, all called {\it Hilbert's reciprocity formula for the cotangent-kernel} (e.g. Erd$\acute{e}$lyi et al.[2], Magnus \& Oberhettinger[3], Muskhelishvili[4]).  The original pair (\ref{26b}) can then be reduced in virtue of (\ref{27}) to perhaps the ultimate form as 
\be  u(\theta)= \hat{H}[v(\phi)]= \frac{1}{2\pi}{\cal P}\int_{-\pi}^\pi v(\phi)\cot \frac{\phi -\theta}{2}d\phi, \quad
     v(\phi) =	\hat{H}^{-1}[u(\theta)]= \frac{-1}{2\pi}{\cal P}\int_{-\pi}^\pi u(\theta)\cot \frac{\theta-\phi}{2}d\theta.\label{28} \nd 
This pair of conjugate equations will be called the {\it circular Hilbert transform}.  
\vskip 0.8mm
\noindent{\bf 6.4. The complementary circular transform.} ~In analogy with the Hilbert transform and its complementary  transform, we can also deduce the transform for the class of function $F(z)$ which is regular in domain ${\cal D}_c^-\{\forall z: ~|z|\geq 1\}$.  For $F(z)$, we simply take (\ref{26a}) or (\ref{28}) with a change in sign of the integral, again implied by (\ref{11}), giving for $F(e^{i\theta})= U(\theta)+iV(\theta)$ the transform equations as
\be  U(\theta)= \check{H}[V(\phi)]= \frac{-1}{2\pi}\int_{-\pi}^\pi V(\phi)\cot \frac{\phi -\theta}{2}d\phi, ~~
     V(\phi) =	\check{H}^{-1}[U(\theta)]= \frac{1}{2\pi}\int_{-\pi}^\pi U(\theta)\cot \frac{\theta-\phi}{2}d\theta,\label{29} \nd 
while the transforms $\hat{H}[\cdot]$ and $\check{H}[\cdot]$ are related exactly like that in (\ref{25b}) with $\hat{H}[\cdot]$ standing for $H[\cdot]$ and $\check{H}[\cdot]$ for $\overline{H}[\cdot]$.  The above pair of relations will be called the {\it complementary circular transform}.
\vskip 0.5mm
\noindent{\it Example 4.} ~As a simple example, we take $v(\phi) = \sin \phi$, then its circular transform, by (\ref{28}), is 
\be  u(\theta)=\frac{1}{2\pi}{\cal P}\int_{-\pi}^\pi \sin(\psi+\theta)\frac{1+\cos\psi}{\sin\psi}d\psi =\cos~\theta, \nonumber \nd
giving $f(e^{i\theta})= u(\theta)+iv(\theta)= e^{i\theta}$, and hence its analytically continued function $f(z)= re^{i\theta}=z$, which is regular in ${\cal D}_c^+$ but is singular at infinity in ${\cal D}_c^-$.  On the other hand, for $V(\phi) = \sin \phi$, we take the complementary circular transform by (\ref{29}), then $U(\theta)=-u(\theta)=-\cos \theta$, giving $F(e^{i\theta})= U(\theta)+iV(\theta)= -e^{-i\theta}$, and therefore $F(z)=-(re^{i\theta})^{-1}= - z^{-1}$, which has a simple zero at $z=\infty$ and is regular inside domain $\overline{{\cal D}_c^-} :(|z|\geq 1)$ as implied by Theorem 10, but is singular at $z=0$ in ${\cal D}_c^+$. 

\vskip 1mm
\noindent{\bf 6.5. The Parseval relations for the generalized Hilbert transforms.}  ~{\it If $u(x)$ and $v(x)$ of (\ref{22}) are both square integrable, and similarly for $U(x), V(x)$ of (\ref{25a}), their complementary counterpart, and further for the pairs $u(\theta), v(\theta)$ of (\ref{28}) and $U(\theta), V(\theta)$ of (\ref{29}), they satisfy the Parseval relations}:
\begin{subequations}\label{30 }
\be  \int_{-\infty}^\infty u^2(x)dx = \int_{-\infty}^\infty v^2(x)dx; && \int_{-\infty}^\infty U^2(x)dx = \int_{-\infty}^\infty V^2(x)dx ;  \label{30a} \\
\int_{-\pi}^\pi u^2(\theta)d\theta = \int_{-\pi}^\pi v^2(\phi)d\phi; &&  \int_{-\pi}^\pi U^2(\theta)d\theta = \int_{-\pi}^\pi V^2(\phi)d\phi. \label{30b} \nd \end{subequations}
\vskip 0.8mm
\noindent{\it Proof:} For the first Parseval relation, we have $u(x)=H[v(t)], ~v(t)=H^{-1}[u(x)]$, then, by (\ref{22}),
\be \int_{-\infty}^\infty u^2(x)dx = \frac{1}{\pi}\int_{-\infty}^\infty u(x)dx\int_{-\infty}^\infty \frac{v(t)dt}{t-x}= \frac{-1}{\pi}\int_{-\infty}^\infty v(t)dt\int_{-\infty}^\infty \frac{u(x)dx}{x-t}= \int_{-\infty}^\infty v^2(t)dt, \nonumber\nd
by interchanging the order of integration.  Similarly, the other Parseval relations can be proved.

Concluding, we note that these various Hilbert transforms all stem from the key relations in (\ref{9}).
\vskip 1mm
\noindent{\bf 7. Behavior of Cauchy function $f(z)$ in the complementary domain $\cal D^-$.}
\vskip 1mm
We have seen exemplified in Example 2-4 that while a Cauchy function $f(z)$ is regular in a closed domain $\overline{{\cal D}^+}$, it invariably has one or more singularities in its complementary domain ${\cal D}^-$.  In general, if $f(z)$ is a Cauchy function, regular inside and on contour $C$, which in this Section will be a unit circle $|t|=1$ for simplicity, then $f(z)$ is implied by Liouville's theorem to possess at least one singularity in $\cal D^-$ outside $C$, including $z=\infty$, unless $f(t) \equiv A (\mbox{const.}) ~\forall t \in C$.  Our primary objective is to determine the exact relationship between the singularities of $f(z)$ in $\cal D^-$ and the values $f(t)$ of $f(z)$ on $C$. 
\vskip 0.5mm
Let us consider here {\it the direct problem}, i.e. with $f(z)$ first prescribed explicitly in terms of all its singularities in $\cal D^-$ outside $C$ in order to examine the corresponding integral formula.  Let one such singularity be located at $z_1, ~|z_1|>1$, which duly induces a corresponding singularity at $z=\infty$, e.g. a pole (or a zero) at $z_1$ inducing a zero (or a pole) of the same order at $z=\infty$; an algebraic or a logarithmic branch point at $z_1$ being associated with the same branch at $z=\infty$, all of which to be accounted for.  These singularities can occur in arbitrary number, of various types, at arbitrary locations $~\forall |z|>1$; their resulting value on $C$, i.e. $f(t) ~\forall |t|=1$ can be deduced at once to be existing and unique.  Of utmost interest is to expound the claim that whatever the system of these singularities of $f(z)$ may be distributed, it invariably results in its corresponding functional $J_{(n,m)}[f(z)]\equiv 0 ~\forall |z|>1$ for various ($m,n$) obeying Theorem (\ref{3b}), (\ref{4b}), (\ref{5b}), and (\ref{6b}), all of which are given by the simple, yet powerful argument of Cauchy's Theorem 1.  To realize this in manifestation, we illustrate it below with a few typical cases of the {\it direct problem}.
\vskip 1mm
\noindent{\it Example 5.}  In (\ref{3a}), given $f(t)=(t-a)^{-1}, ~(|a|>1)$, which is a simple pole situated outside $C ~(|t| =1)$ and a simple zero at $z=\infty$, so we have, by (\ref{3a}), $J[f(z)]= (z-a)^{-1}=f(z)$ for $|z|\leq 1$, which is regular for $|z|\leq 1$ (since $|a|>1$), whereas for $|z|>1$, 
\be 2\pi i J[f(z)] &=&\oint_C g(t,z)dt = I_C =0 \qquad (g(t,z)=(t-a)^{-1}(t-z)^{-1}, ~|a|>1, |z|>1) \nonumber\\
   &=& I_{\infty}-I_S = \left\{\oint_{C_\infty}-\oint_{C_S}\right\}g(t,z)dt = \frac{1}{z-a}\oint_{C_S} (\frac{1}{t-a}-\frac{1}{t-z})dt = 0, \notag\nd
which results as follows.  The three integrals, $I_C$ on contour $C (|t|=1)$, $I_{\infty}$ on $C_\infty$ encircling $t=\infty$, and $I_S$ enclosing both poles of $g(t,z)$ (all in the positive sense), are related by (\ref{1}) as $I_{\infty}= I_C+I_S $ since $g(t,z)$ is regular in the domain bounded by $C_\infty, C$ and $C_S$.  Separately, $I_C=0$ by (\ref{3b}) (g(t,z) being regular $\forall~|t|\leq 1, |a|>1, |z|>1$), $I_{\infty}=0$ since $g(t,z)=(t^{-2}+O(|t|^{-3}))$ with residue $res.=0$ at $z=\infty$, and $I_S=0$ since the residues of its integrand at $t=a ~(res.=1)$ and at $t=z ~(res.=-1)$ cancel.  Whence $J[f(z)]=0$ for $|z|>1$ is shown both by applying (\ref{3b}) and alternatively by direct integration. 
\vskip 1mm
\noindent{\it Example 6.}  As a versatile variation, let the complement function in (\ref{7a}) be $F(t)=t^{-n} ~(n= 1, 2, \cdots)$ which has a pole of order $n$ at $t=0$ and is regular for $|t|\geq 1$, then for $|z|\geq 1$ we have $J^-[F(z)]= I_{C^-}=(2\pi i)^{-1} \oint_{C^-} t^{-n}/(t-z)dt = z^{-n}$ by (\ref{7a}), or by direct evaluation, $I_{\infty}=I_S + I_C =0 $ (due to its zero residue at $t=\infty$), hence $I_{C^-}=-I_C=I_S = z^{-n}$ (the residue at $t=z, |z|>1$).  For $|z|<1$, 
\be 2\pi i J^-[F(z)] &=& \oint_{C^-} g(t,z)dt = I_{C^-}=0 \qquad (g(t,z)=t^{-n}(t-z)^{-1},~|z|<1), \nonumber\\
   &=& I_{C^-}=-I_C=-I_\infty = 0, \nonumber\nd  
where the first result for $I_{C^-}=0 ~(|z|<1)$ is by (\ref{7b}) of Theorem 6, whereas the second results from direct integration by deforming the contour $C$ to $C_{\infty}$, between which $g(t,z)$ is regular and $g(t,z)=O(|t|^{-(n+1)})$ with zero residue at $t=\infty$, hence the result.
\vskip 1mm
\noindent{\it Example 7.}  In (\ref{3a}), given $f(t)=(t-1/b)^{-1/2} ~(0<b<1)$, made single-valued on a two-sheet Riemann surface cut along the real $t$-axis from branch point at $t=1/b>1$ to $t=+\infty$, hence by (\ref{3a}), $J[f(z)]=(z-1/b)^{-1/2}$ which is regular for $|z|\leq 1$ ($1/b>1$), whereas for $|z|>1$,    
\be 2\pi i J[f(z)] &=&\oint_C g(t,z)dt = I_C =0 \qquad\qquad\quad (g(t,z)=(t-1/b)^{-1/2}(t-z)^{-1}, ~|z|>1) \nonumber\\
   &=& \oint_{C_I} g(\frac{1}{\xi}, z)\frac{d\xi}{\xi^2} = \frac{2\sqrt{b}}{z}\left\{\int_0^b \frac{d\xi}{\sqrt{\xi(b-\xi)}(1/z -\xi)}- \frac{\pi z}{\sqrt{1-bz}}\right\}=0 \quad (1<|z|<1/b), \nonumber\nd
where $I_C=0$ by (\ref{3b}), $\xi=1/t$ is the inverse mapping, with contour $C_I ~(|\xi|=1)$ in the positive sense, the line integral from $\xi=0$ to $b$ comes from the contour integral around the cut within $C_I$, whilst the last term in the bracket comes with the residue at $\xi=1/z$ (located within $C_I$).  We note that the value $I_C=0$ (or equivalently, $\oint_{C_\infty}g(t,z)dt=0$) is necessary and sufficient to have the line integral determined as shown above.  A similar result can be attained when point $z$ falls on the branch cut.
\vskip 1mm
\noindent{\it Example 8.} ~Finally, let us consider the special case with $f(t)\equiv 1$ on $|t|=1$, for which we have
\be J[f(z)]&=& \frac{1}{2\pi i} \oint_C \frac{dt}{t-z}= I_C=1  \qquad\qquad\quad  (|z|\leq 1), \nonumber \\
	&=& I_{\infty}-I_S = \frac{1}{2\pi i} \oint_{C_\infty}(\frac{1}{t}+O(|t|^{-2}))dt-\frac{1}{2\pi i} \oint_{C_S}\frac{dt}{t-z}= 1-1=0 \quad (|z|>1).  \nonumber \nd
We note that in this case, even with the functional $J[f(z)]\equiv 1 ~(|z|\leq 1)$, $J[f(z)]$ still jumps down to $J[f(z)]\equiv 0 ~(|z|> 1)$ in fulfilling (\ref{3b}) of Theorem 2.
\vskip 1mm
Therefore we can claim that only when $f(z)= 0$ on contour $C$ do we achieve the unique particular result of {\it having functional $J[f(z)]\equiv 0$ uniformly continuous throughout the entire $z$-plane, while Theorem (\ref{3b}) invoking $J[f(z)]\equiv 0$ for $z\in ~{\cal D}^-$ is universally fulfilled}.  We thus recognize the sharp contrast between any function $f(z)$ which is regular $\forall z \in {\cal D}^+$ and its contour-integral functionals $J_{(n,m)}[f(z)]$ in their characteristic behavior regarding their continuity, convergence, and singularity distributions over the $z$-plane; and so does this hold for complement function $F(z)$.  Conceptually, it serves no further purpose to pursue the functional $J[f(z)]$ as a function, $H(z)$ say, for it is no more than $H(z)=f(z)$ being regular in domain $\overline{{\cal D}^+}$ and $H(z)\equiv 0$ in domain ${\cal D}^-$.  $H(z)=J[f(z)]$ is simply not an analytic function; it is neither continuous nor differentiable in a neighborhood striding across contour $C$.

\vskip 1mm
\noindent{\bf 8. The Plemelj formulas.}  ~We now consider another general class of line integrals of the form
\be   f(z)= \frac{1}{2\pi i}\int_L \frac{g(t)}{t-z}dt \qquad (z \notin L),\label{31}\nd
where $L$ is a regular Jordan arc, assumed finite (or infinite) in length, with end-points at $t=a$ and $t=b$ and without double point, the integration variable $t$ moves from $a$ to $b$~(see Fig. 2), and $g(t)$ is assumed regular in a neighborhood of the entire line $L$.  Evidently, $f(z)$ is a regular function $\forall z\notin L$ and has a simple zero at infinity.  It is also evident that $f(z)$ possesses derivatives to all orders, given by
\be   f^{(n)}(z) = \frac{n!}{2\pi i}\int_L \frac{g(t)dt}{(t-z)^{n+1}} \qquad (z \notin L, ~n=1, 2, \cdots). \label{32}\nd

\begin{figure}[ht]
\begin{center}
\includegraphics[scale=0.64, trim=16mm 16mm 12mm 15mm]{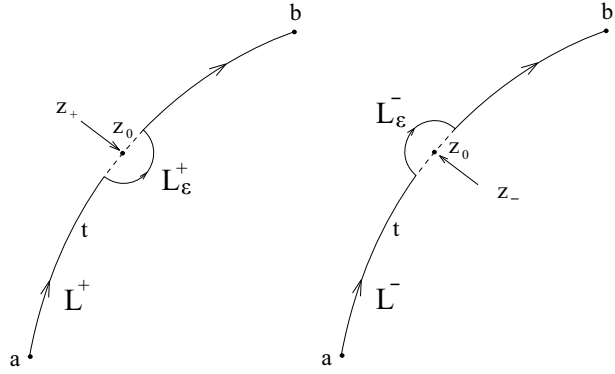}
\caption{\small  A Plemelj integral $f(z)=(2\pi i)^{-1}\int_L g(t)dt/(t-z)$ along an open path $L$ (from $t=a$ to $t=b$) has its path $L$  deformed into $L^{\pm}$ with only an indented small semicircle $L_{\epsilon}^{\pm}$ ($|t-z_0|=\epsilon$), on the ${\mp}$ive side of $L$, centered at $z_0\in L$ so as to let a point $z$ on the ${\pm}$ive side (left or right side) of $L$ tend to $z_0\in L$ without crossing $L^{\pm}$.}
\end{center}
\end{figure}
\vskip -3.8mm

In the limit as point $z$ tends to a point $t=z_0$ on $L$ from the left (+ive), or from the right (-ive) side of $L$, we indent $L$ into $L^\pm=L^{\pm}_\epsilon+L_{p}$ where $L^{\pm}_\epsilon$ is a semicircle of radius $|t-z_0|=\epsilon$ onto the ({$\mp$}ive) side, leaving $L_{p}=L^{\pm}-L^{\pm}_\epsilon$ intact, so that point $z$ reaches $z_0\in L$ without crossing $L^{\pm}$, while $f(z)$ tends to its limit $f^+(z)$, or $f^-(z)$, respectively, yet undetermined.   Carrying out the integration over $L^{\pm}$ paths in a way similar to that for $C^{\pm}$ in \S 3 resulting in (\ref{8a})-(\ref{8b}), we obtain Plemelj's formula: 
\be  f^{\pm}(z)= \pm \frac{1}{2}g(z) +\frac{1}{2\pi i}{\cal P}\int_{L} \frac{g(t)}{t-z}dt \qquad (z \in L - \mbox{Plemelj's formula}),\label{33}\nd
where the suffix of $z_0$ is omitted, the sign ${\cal P}$ (often omitted) signifies its Cauchy principal value,
\be  {\cal P}\int_{L} \frac{g(t)}{t-z}dt = \lim_{\epsilon\rightarrow 0}\left\{\int_a^{z-\epsilon} +\int_{z+\epsilon}^b\right\}\frac{g(t)}{t-z}dt \qquad (z \in L),\notag\nd
where both $z-\epsilon$ and $z+\epsilon$ lie on the regular arc L.   ~From (\ref{33}) we deduce
\begin{subequations}\label{34 } 
\be  f^{+}(z)- f^{-}(z) &=& g(z)  \qquad\qquad \qquad\quad(z \in L), \label{34a} \\
     f^{+}(z)+ f^{-}(z) &=& \frac{1}{\pi i}{\cal P}\int_{L} \frac{g(t)}{t-z}dt ~\qquad (z\in L).\label{34b} \nd \end{subequations}
Formulas (\ref{33}) and (39a,b) are called {\it Plemelj's formulas}.  Finally, substituting (\ref{34a}) in (\ref{31}) yields 
\be   f(z)= \frac{1}{2\pi i}\int_{L} \frac{g(t)}{t-z}dt =\frac{1}{2\pi i}\int_{L} \frac{f^{+}(t)- f^{-}(t)}{t-z}dt ~,\label{35} \nd   
which shows that $f(z)$ is determined by its jump [$f^{+}(t)- f^{-}(t)$] across the line $L$ for all $z$ in the plane, including the $z$'s on $L$, on which (\ref{35}) reduces to an identity by virtue of (\ref{34a})-(\ref{34b}).  

In applying Plemelj's formulas to integrals involving two Cauchy kernels, as exemplified in Example 1, it is essential for analysis in application to have {\it the Poincar$\acute{e}$-Bertrand formula}: 
\be  \int_L\frac{dt'}{t'-z_0}\int_L\frac{f(t,t')}{t-t'}dt = \int_L dt\int_L\frac{f(t,t')dt'}{(t'-z_0)(t-t')} - \pi^2 f(z_0,z_0) \qquad  (x_0 \in L), \label{36} \nd
in which each integral assumes its own principal value, here with the symbol $\cal P$ omitted as understood by convention.  For its proof we can let a point $z\notin L$ tend to a point $z_0\in L$, while applying Plemelj's formula twice as needed, as can be found shown in the literature (e.g. Muskhelishvili[4]).
\vskip 0.8mm
We further remark that proofs for (\ref{33}) with $g(t)$ more restricted can be found (e.g. Muskhelishvili[4]) with $g(t)$ assumed to satisfy the H$\ddot{o}$lder condition, i.e. for $z, t$ on $L$, $|g(t)-g(z)|<M|t-z|^{\mu}$ for some positive constants $M$ and $\mu$, $\mu$ being the H$\ddot{o}$lder index.  The theorem was proved by Plemelj (1908) for $0<\mu <1$.  In these aspects, Plemelj's formulas have powerful applications to various two-dimensional mathematical physics, including problems with mixed boundary conditions on harmonic and biharmonic functions, the Riemann-Hilbert problems and other types of problems encountered in fluid and solid mechanics, elasticity, physics, fields of engineering and applied mathematics.

\vskip 1mm
\noindent{\bf 9. Application and physical significance.} ~Applications of the formulas obtained in this study can be made to various scientific fields as well as for further mathematical advances.  Here we first select a Riemann-Hilbert problem of aerodynamic wing theory for an exact solution to a two-dimensional flat plate airfoil.  It involves dealing with a pair of conjugate integrals, one of which is about a {\it finite Hilbert transform and its inversion}, whilst the other a Plemelj's integral over an infinite line.

Thus we consider the 2-D irrotational flow of an incompressible and inviscid fluid past a flat plate airfoil held fixed along $-1\leq x\leq 1, y=0$ in an inertial frame of reference at an incidence angle $\alpha$ with respect to a uniform free stream of velocity $\mb{U}$.  Denoting the fluid velocity by $(U \cos \alpha +u, U \sin \alpha +v), ~(u,v)$ being the perturbation velocity, we have the basic equations[9] as 
\begin{subequations}\label{37 }
\be   u_x +v_y =0  & &(\mbox{incompressibility})  \label{37a} \\
   u_y -v_x =0     & &(\mbox{irrotationality})  \label{37b} \\
  v= -U\sin \alpha \quad (-1\leq x\leq 1, y=\pm 0);  & &u^2+v^2 \rightarrow 0 \quad (\mbox{as} ~x^2+y^2\rightarrow \infty),  \label{37c} \\
 \frac{p}{\rho} +\frac{1}{2}\{(U\cos\alpha +u)^2+(U\sin\alpha +v)^2\}&=& \frac{1}{2}U^2,\label{37d} \nd \end{subequations}
where $p$ is the fluid pressure, gauged to zero at infinity, and $\rho$ the fluid density.  Here (42a,b) are the 2-D components of $\nabla\cdot\mb{u}=0$ and $\nabla\times\mb{u}=0$, respectively, ~$\mb{u}$ being the vector $(u,v,0)$.  With (\ref{37c}) providing the boundary conditions, a solution to $(u,v)$ can be found from (42a-c), and (\ref{37d}) then gives pressure $p$ by this Bernoulli equation.  This completes the mathematical formulation of the problem. 

Noting that ($u, -v$) satisfy the Cauchy-Riemann equations (42a,b), the {\it complex velocity} $w=u-iv$ is then an analytical function of $z=x+iy$.  Since $v(x,y)$ is prescribed in (\ref{37c}) as being even in $y$, then $u(x,y)$ by (42a,b) is odd in $y$.  Since $w(z)$ is analytic and regular in the open domain of the flow field, it follows that $u(x,0)=0$ for $|x|>1$ since $u(x,y)$ is there continuous and odd in $y$.  As $v(x,0)$ is prescribed in $|x|<1$ and $u(x,0)=0$ for $|x|>1$, this is a Riemann-Hilbert problem.  It has a complementary solution, namely $w_c=u_c-iv_c=i/H(z), H(z)=\sqrt{z^2-1}$, for on $y=0$, $H^{\pm}(x)=\pm i\sqrt{1-x^2} ~(|x|<1)$ and $H^{\pm}(x)=\mbox{sgn} ~x ~\sqrt{x^2-1} ~(|x|>1)$, so that $v_c =0$ for $|x|<1$, $u_c= 0$ for $|x|>1$, and $|w_c|\rightarrow 0$ as $|z|\rightarrow \infty$.  Now introducing $f(z)=w(z)H(z)$ gives
\be  f^+(x)-  f^-(x) &=& 2\sqrt{1-x^2}~v(x) \qquad\quad~ (|x|<1, ~v(x)=-U\sin\alpha ), \notag\\
			&=& 0 \qquad\qquad\qquad\qquad\quad ~(|x|>1). \notag\nd
Hence, by Plemelj's formula (\ref{35}) (here with the path $L$ spanning the entire $x$-axis), we have
\be  w(z)= \frac{1}{\pi i \sqrt{z^2 -1}}\int_{-1}^1 \frac{\sqrt{1-t^2}}{t-z}~v(t)~dt + \frac{iB}{\sqrt{z^2 -1}}, \notag\nd
B being an arbitrary real constant.  This is the only correct form for a solution to $w(z)$ if 
$|w(z)|=O(|z|^{-1})$ as $z\rightarrow \infty$ and $w(z)$ be integrable at the plate.  Finally, $B$ is determined by an additional {\it physical condition}, known as {\it Kutta's condition}, requiring[9] that $w(z)$ be regular in a neighborhood of the trailing edge at $z=1$.  Hence expanding this $w(z)$ about $z=1$ yields the unique exact solution as 
\be  w(z) = -\frac{1}{\pi i} \sqrt{\frac{z-1}{z+1}}\int_{-1}^1 \sqrt{\frac{1+t}{1-t}}~\frac{v(t)}{t-z}~dt \qquad (0\leq |z|<\infty), \label{38} \nd
valid for arbitrary $v(x)$.  For the flat plate, $v(x)=-U\sin\alpha$ by (\ref{37c}), (\ref{38}) reduces by using (\ref{33}) to
\be  u^{\pm}(x)-i v^{\pm}(x)= U\sin \alpha ~\left(\pm \sqrt{\frac{1-x}{1+x}}+i\right)  \qquad (-1<x\leq 1), \label{39} \nd
exhibiting that $u$ has an equal but opposite jump across the plate by a distribution having a square root singularity at the leading edge at $z=-1$ and vanishing at the trailing edge at $z=1$.  This jump distribution of $u$ results in the so-called {\it circulation}, $\Gamma$, around the plate (see, e.g. von K$\acute{a}$rm$\acute{a}$m \& Burgers[9]), given by the contour integral of $u$ clockwise (by convention) around the airfoil, 
\be  \Gamma = \oint u(x, \pm 0)dx = 2\pi U\sin \alpha.  \notag \nd

Finally, expressed in three-dimensional vectors, $\mb{U}=(U\cos\alpha, U\sin\alpha, 0), \mb{\Gamma}=(0, 0, -\Gamma)$ (by the right-hand rule with the contour integral for $\Gamma$), and $\mb{L}=(L_1, L_2, 0)$ for the lift vector acting on the airfoil, we have lift $\mb{L}$ given by the Kutta-Joukowski theorem[9] in vector cross product of $\mb{U}\times\mb{\Gamma}$ as 
\be   \mb{L}=\rho \mb{U}\times\mb{\Gamma}, \qquad  L=|\mb{L}|= 2\pi\rho U^2 \sin \alpha,  \label{40} \nd
by which the lift $\mb L$ acts perpendicular to the free stream velocity $\mb U$, pointing upward if positive.
\vskip 1mm
Regarding this problem and its solution, there are several issues worthy of expository discussion.
\vskip 0.5mm
\noindent{\bf 9.1. Physical significance and advances in applications.}  Mathematically, the above solution to the airfoil problem as formulated is {\it exact}.  It may have value in showing a methodology useful for achieving exact solutions and in serving as a standard reference for assessing approximate approaches such as by linear theory for small incidence angles.  Physically, however, the square root singularity of velocity at the leading edge (associated with an even worse singular suction in pressure $p$ by (\ref{37d})) should raise serious questions concerning not only for engineering applications but further for experimental verification of the range of validity, even for very small incidence angles.  Such concerns have actually stimulated innovative theoretical developments by deftly superposing a distribution of flow-mass sources along the plate so as to obtain, again in exact form, real airfoil profiles enclosing the lifting flat plate and its singularity (lying inside the profile now having a round nose at the leading edge and a cusped trailing edge), with results that can be subjected to specific engineering design and wind tunnel tests for validation and adaptation by the industry.  In return, the result is gratifying that some airfoils with aptly designed round noses have been found experimentally capable of sustaining the low suction pressure for incidence angle $\alpha$ up to around 18 degrees before the airfoils stall, so to speak, with flow separation.  Another fruitful reward is the experimental support to Kutta-Joukowski's theorem that the inviscid lift vector $\mb L$ is well predicted both in magnitude and direction, with the small viscous effects accounted for.  This approach to determine accurate solutions to problems of fluid flow past bodies of finite volume by placing flow singularities at an ultimate focal point or plane (e.g. at the center of a circle or sphere or at the focal ellipse of a tri-axial ellipsoid) is now classical for mechanics and electrodynamics, followed by more advanced methods for aerodynamics of thin airfoils (e.g Lighthill [10]), for naval hydrodynamics of double-body for ship hulls (e.g. Wu \& Chwang [11]), for biharmonics of triaxial ellipsoids and in other cases.  In this respect, having the exact solution can provide a concrete foundation in general for further sound overall development.

Returning to the mathematics, the infinite suction pressure acting at the pointed leading edge can indeed be integrated with rigor to produce a so-called {\it finite leading edge suction} $\mb{S}$ pulling the airfoil forward along the flat plate just so exactly as to make the resultant lift $\mb L$ (as the vector sum of $\mb{S}$ and the pressure integral acting normal to the plate) to act normal to the free stream velocity $\mb U$, as predicted by Kutta-Joukowski theorem.
\vskip 0.5mm
\noindent{\bf 9.2. Finite Hilbert transform and its inversion.} ~When Cauchy integrals enter the analysis for a physical problem, the primary query would be on the significance of the Cauchy kernel $(t - z)^{-1}$.  In mathematical physics, Cauchy integrals are closely related to the potentials of single and double layers distributed along a contour $C$ or an open arc $L$, as we can illustrate next.  For this airfoil problem, there are actually alternative approaches for the solution.  We have elected above to regard it as a Riemann-Hilbert problem.  But we can also represent the flat plate in the free stream by a distribution of flow singularities fixed to the plate in the so-called {\it singularity method}.  For 2-D flows satisfying (42a,b), it is well known that the elementary flow singularities are a point source of strength $Q$ and a point vortex of strength $\Gamma$ held at the origin, say, giving their own velocities at a field point $z$ as
\be  w(z)=u(x,y)-iv(x,y)= \frac{Q+i\Gamma}{2\pi z} \quad\longrightarrow\quad \overline{w(z)}=u+iv=\frac{Q-i\Gamma}{2\pi r}e^{i\theta}, \label{41}\nd
which, now expressed in the polar coordinates, $z=x+iy=r e^{i\theta}$, shows that source $Q$ has only an outward radial velocity component $u_r=Q/2\pi r$ ($\arg \overline{w}=\theta$), and vortex $\Gamma$ has only a clockwise circumferential velocity component $u_{\theta}=-\Gamma/2\pi r$ ($\arg \overline{w}=\theta-\pi/2$).  In terms of these base singularities, we can construct a surface distribution of sources of density $q$ and of a vortex sheet of density $\gamma$ per unit length along a regular arc L, generating a complex velocity field $w(z)$ as
\be  w(z) = \frac{1}{2\pi}\int_{L} \frac{q(t)+i\gamma (t)}{z-t}dt, \label{42} \nd
here with the Cauchy kernel physically signified.  For the flat plate airfoil held fixed in a free stream, we need only a surface distribution of vortex sheet along the plate, giving its complex velocity as 
\be  w(z) = \frac{1}{2\pi i}\int_{-1}^1 \frac{\gamma (t)}{t-z}dt. \label{43} \nd
From this we have, by Plemelj's formula (\ref{33}), that on the $\pm$sides of the plate,
\begin{subequations}\label{44 }
\be  w^{\pm}(x)= u^{\pm}(x)-i v^{\pm}(x)&=& \pm\frac{1}{2}\gamma (x)+ \frac{1}{2\pi i}\int_{-1}^1 \frac{\gamma (t)}{t-x}dt \qquad (|x|<1), \nonumber \\
  \longrightarrow\quad  u^+(x)-u^-(x) &=& \gamma (x) \qquad\qquad\qquad\qquad\qquad\quad ~~(|x|<1), \label{44a} \\
     v^+(x)=v^-(x) &=& \frac{1}{2\pi }\int_{-1}^1 \frac{\gamma (t)}{t-x}dt\equiv G[\gamma(t)]. \qquad ~(|x|<1).\label{44b} \nd \end{subequations}
With the boundary conditions (\ref{37c}) prescribing $v(x)$, (\ref{44b}) actually is a singular integral equation for the vorticity distribution $\gamma (x)$, which mathematically can also be regarded as a {\it finite Hilbert transform}, with the integral operator $G$ denoting the transform.  Then its inversion can be given by (\ref{44a}), with $u^+ -u^-= \gamma (x)$ deduced from (\ref{38}) for arbitrary $v^\pm (x)=v(x)$, yielding the unique solution as 
\begin{subequations}\label{45 }
\be  v(x) &=& \frac{1}{2\pi }\int_{-1}^1 \frac{\gamma (t)}{t-x}dt\equiv G[\gamma(t)] \quad\longrightarrow\quad \gamma (x)=G^{-1}[v(t)]. \qquad (|x|<1), \label{45a}\\
   \gamma (x) &=& -~\frac{2}{\pi}\sqrt{\frac{1-x}{1+x}}\int_{-1}^1\sqrt{\frac{1+t'}{1-t'}}\frac{v(t')}{t'-x}dt'\equiv G^{-1}[v(t')]\qquad\qquad~(|x|<1). \label{45b}\nd \end{subequations}
where the integral operator $G$ defines the {\it finite Hilbert transform} and the operator $G^{-1}$ its inversion, signifying $GG^{-1}=G^{-1}G=1$ (the unity operator).  This is equivalent to substituting (\ref{45b}) into the integral equation (\ref{45a}) for a final check, which can be shown for arbitrary $v(x)$ by interchanging the order of integration by applying the Poincar$\acute{e}$-Bertrand formula (\ref{36}) as shown in Example 1. 
\vskip 0.5mm
\noindent{\bf 9.3. Generalizations for application.} ~
\vskip 0.5mm
This primary application could provide a sound basis for generalizations as often pursued in fluid mathematics.  We will only describe briefly here an extension of the stationary flat plate as a base lifting surface to develop a fully nonlinear theory for a two-dimensional flexible wing moving with arbitrary unsteady variations in wing profile and along arbitrary trajectory for modeling bird/insect flight and fish swimming.

Thus, we consider the irrotational flow of an incompressible and inviscid fluid produced by a two-dimensional flexible lifting surface $S_b(t)$ of negligible thickness, moving with time $t$ through the fluid in arbitrary manner.  Its motion is described by using a hybrid Lagrangian-Eulerian system with the Lagrangian body coordinates ($\xi,\eta$) to identify a point $X(\xi, t), Y(\xi,t)$ at time $t$ on the wing surface $S_b(t)$ and on the vortex sheet $S_w(t)$ shed from the wing, both of which can be prescribed by a complex coordinate $z=x+iy$ (for the Euler description) fixed in an inertial frame of reference, and with $z = Z(\xi, t)$ prescribed for the time-dependent body-wake motion function, parametrically in $\xi$ as 
\be  Z(\xi, t)= X (\xi, t) + i Y(\xi, t) \quad\mbox{on}\quad  S_b(t):(-1<\xi<1)+ S_w(t):(1<\xi<\xi_{m}) , \label{46}\nd 
with the leading and trailing edges of the wing at~$\xi =-1$ and $\xi = 1$, respectively, while the vortex sheet is shed again smoothly from the wing trailing edge under the Kutta condition to form a prolonging wake $S_w(t)~(1<\xi\leq\xi_{m}$ with $Z(\xi_{m}, t)$ charting the position of the starting vortex shed at $t=0$ to reach $\xi_m =\xi_{m}(t)$ at time $t$.  The problem is formulated with a boundary-value requiring the flow velocity normal to $S_b(t)$ equal to that of $S_b(t)$ itself and with the initial value when the motion starts at time $t=0$ in an unbounded fluid at rest in the inertial frame of reference and with $S_b(t=0)$ in a stretched-straight shape such that $Z(\xi, 0)=\xi ~~(-1<\xi<1, ~\eta =0)$.  For $t>0$, the point $\xi$ on $S_b(t)$ moves with a {\it prescribed} body motion function $Z(\xi,t)$ and a {\it prescribed} complex velocity $W(\xi, t) = U-iV$, 
\be  W(\xi, t) = U-iV = \partial\overline{Z}/\partial t = X_{t}-i Y_{t} \qquad (|\xi|<1, ~t\geq 0; ~\overline{Z}=X-iY), \label{47}\nd
which has a tangential component, $U_{s}(\xi,t)$, and a normal component, $U_{n}(\xi,t)$, given by
\be  W\partial Z/\partial \xi =(X_{\xi}X_t +Y_{\xi}Y_t)-i(X_{\xi}Y_t -Y_{\xi}X_t)=U_{s}-iU_{n} \label{48}\nd
on $S(t)=S_b(t)+S_w(t)$, the flexible $S_b(t)$ being assumed inextensible ($|\partial Z/\partial \xi|=1$).  Thus, the normal velocity, $U_{n}(\xi,t)$, is prescribed for $|\xi|<1$ while the wake vortex, once shed, is conserved as free vortex in motion, and the problem is to determine the vorticity $\gamma (\xi,t)$ over $S(t)=S_b(t)+S_w(t)$.  For the solution, we refer to the studies by Wu[12] for the details.

\vskip 0.8mm
\noindent{\bf 10. Discussion and conclusion.}
\vskip 0.8mm
The primary objective of the present study is first to extend the coverage of Cauchy's integral formula (1a,b) to include the contour $C$ of the integral so as to render it valid for the entire $z$-plane.  With Cauchy's function $f(z)$ assumed $C^n~\forall z\in {\cal D}^+$ within $C$ and about $C$, $f(z)$ and all its derivatives $f^{(n)}(z)$ are proved to be uniformly continuous in the closed domain $\overline{{\cal D}^+}=[{\cal D}^+ +C]$.  Under the same assumption, Cauchy's integral formulas (for $n=0,1,\cdots$) $\forall z\in {\cal D}^+$ (or $\forall z\in {\cal D}^-$) are proved uniformly convergent in closed domain $\forall z\in\overline{{\cal D}^+}$ (or $\forall z\in\overline{{\cal D}^-}=[{\cal D}^- +C]$).  From these fundamental discoveries there follow findings of various integral properties of $f^{(n)}(z)$ and $J_n[f(z)]$.  A new complement function $F(z)$ is introduced to be $C^n~\forall z\in {\cal D}^-$ outside $C$ and about $C$, and shown to share all the claims for $f(z)$ in complete analogy.  These new results have provided a simple and sound base to derive the generalized Hilbert transforms in various domains of different geometry, and to explore jointly the roles of Plemelj's formulas in application to engineering science, mathematical physics, and applied mathematics.  In conclusion, there are nevertheless several vital issues of great significance calling for expository discussion.
\vskip 0.6mm
\noindent{\bf 10.1. Overall behavior of the Cauchy function in the entire $z$-plane.}  ~The foregoing deliberation on the general properties of Cauchy function $f(z)$ over the entire $z$-plane is comprehensively expounded by considering the {\it direct problem}, i.e. with $f(z)$ first prescribed by assigning an arbitrary distribution of all its zeros and singularities in domain ${\cal D}^-$ outside $C$, including $z=\infty$, as seen exemplified in Examples 5-8 of \S 7.  With $f(z)$ thus prescribed in explicit functional expressions, it is obvious that not only the values of $f(t)$ are all known $\forall t\in C$, but also the value of $f(z)$ is uniquely given $\forall ~z \in {\cal D}^+$ without having to use the integral formula for evaluation; and even much further, the singularities of $f(z) ~\forall z\in {\cal D}^-$ are all completely revealed, as exemplified in Examples 5-8.  In general practice, however, resolving mathematical problems encountered in engineering science, mathematical physics, and applied mathematics, formulated as certain initial-boundary value problems or with some integro-differential equations for numerical computation, resort is often to take certain suitable approach, e.g. a perturbation expansion scheme with unknown coefficients calculated in numerics.  In such countless cases for 2-D studies, the primary variable could be a function $f(z)$ of a complex variable $z$, regular in the problem domain ${\cal D^+}$ and only known numerically $\forall z \in {\cal D^+}$ and on its boundary contour $C$.  Then Cauchy's integral formula is capable of providing such values as $J[f(z)]=f(z)~\forall z\in \overline{{\cal D^+}}$, $J[f(z)]\equiv 0~\forall z\in {\cal D^-}$ outside $C$, and their derivatives by our new Theorems 6-8, yet providing no clue at all for $f(z)$ in ${\cal D^-}$.  It is to this end that we have need to proceed onto the following issues. 
\vskip 0.8mm
\noindent{\bf 10.2.  Relationship between the Cauchy integral and Plemelj's formulas.} ~In sharp contrast to the overall behavior of Cauchy function $f(z)$ and its integral $J[f(z)] ~\forall z ~(0\leq |z|<\infty)$ as just expounded above, a Plemelj integral, $f(z)=(2\pi i)^{-1}\int_L g(t)dt/(t-z)$ along a regular Jordan arc $L$ without a double point (not closed as a contour) is actually an analytic function, regular $\forall z\notin L$ and has a simple zero at $z=\infty$.  Only in the limit of a point $z\notin L$ tending to a point $z_0\in L$ from the opposite sides of $L$ does the function $f(z)$ endure an equal and opposite jump (i.e. equal to $\pm g(z_0)/2$ on the $\pm$-sides) and at the same time with $f(z)\rightarrow f^{\pm}(z_0)$ in the limit as shown in Plemelj's formula (\ref{33}).  However, once the two ends of $L$ coincide to form a closed simple contour $C$ enclosing an open domain $\cal D^+$ bounded by $C$ and excluding an open domain $\cal D^-$ outside $C$, then immediately the Plemelj integral becomes a Cauchy integral, $J[f(z)]$, conjointly associated with the stark change in $f^{+}(z_0)$ to become the prescribed $f(z_0) \forall z_0\in C$ together with $f^{-}(z_0)\equiv 0$, by our Theorem 2.  This conspicuous change in value of the integral when an open integration path $L$ is changed into, or from, a closed contour $C$ is of course a natural consequence to the mathematics in exact rigor, yet still seems greatly worthy of calling for special attention.  We reiterate that the functional $J[f(z)]$ is simply not an analytic function; it is neither continuous nor differentiable in a neighborhood striding across contour $C$.
\vskip 0.8mm
\noindent{\bf 10.3. A conjecture on an unsolved inverse problem.} ~Finally, we conclude the present study with high expectation in bringing forth an inverse problem of great significance as follows.  This is concerned with the generalized Cauchy function $f(z)$ associated with Cauchy's integral formula (\ref{3a})-(\ref{3b}) involving integral $J[f(z)]$ over a contour $C$ in special regard to the relationship between all its singularities $\forall z \in {\cal D^-}$outside $C$ and the values of $f(t)~\forall t\in C$ which are known only {\it numerically}, or in terms of a series with known numerical coefficients, but NOT in any closed functional expression as exemplified in Examples 5-8.  The relationship is obviously known for the direct problem as explained, but seems to require resolution for the inverse problem.
\vskip 0.8mm
\noindent{\bf The inverse problem.} ~{\it The inverse problem is to have only the numerical data given for $f(t)~\forall t \in C$ for a function $f(z)$ being regular inside contour $C$ and using them to determine all the exact singularity distribution of $f(z)~\forall z\in \cal D^-$ outside $C$ in a closed analytical form, whatever the singularity distribution.}  
\vskip 0.8mm
\noindent{\bf The conjecture.} ~{\it We conjecture that a solution to this inverse problem exists, which may not be unique.}  
\vskip 0.5mm
It is hoped that the direct problem delineated in \S 10.1 and exemplified in \S 7 with Example 5-8 may cast light on the course to its resolution.  Having a general methodology to the solution of this inverse problem is of vital importance, for in studies of physical phenomena, solutions are usually found in numerics, to have solution to this inverse problem is essential to gaining in-depth comprehension of the phenomenon in pursuit, such as feasible instabilities and bifurcations of the solution found in existence.
\vskip 1mm
\noindent{\bf Acknowledgment.} ~I wish to thank Prof. Joe Keller, Prof. John C.K. Chu, Prof. Michael Weinstein, and Prof. Lu Ting for interesting discussions, and especially Prof. Jin Zhang of Hong Kong University for careful reading of the text and analysis.  I am most appreciative for the gracious encouragement from Dr. Chinhua S. Wu and the American-Chinese Scholarship Foundation.

\vskip 1.8mm
\noindent{\bf References.}
\vskip 0.8mm
\noindent [1] Courant, R. \& John F. Introduction to Calculus and Analysis.  Interscience Publ. ~(1974). 
\vskip 0.5mm
\noindent [2] Erd$\acute{e}$lyi, A. (Editor) Tables of Integral Transforms, Vol. 2 (Bateman Manuscript Project).   MaGraw-\\ \indent Hill Inc. ~(1954).
\vskip 0.5mm
\noindent [3] Magnus, W. \& Oberhettinger, F. Special Functions of Mathematical Physics. Chelsea Pub. ~(1949).
\vskip 0.5mm
\noindent [4] Muskhelishvili, N.I. Singular Integral Equation.  Noordhoff N.V. ~(1953).
\vskip 0.5mm
\noindent [5] Stokes, G.G. Mathematical and Physical Papers, Vol. 1, 5. Cambridge U. Press~(1880).
\vskip 0.5mm
\noindent [6] Titchmarsh, E.C. The Theory of Functions.  Oxford U. Press ~(1949).
\vskip 0.5mm
\noindent [7] Titchmarsh, E.C. The Theory of Fourier Integrals.  Oxford U. Press ~(1948).
\vskip 0.5mm
\noindent [8] Wu, Th.Y., Kao, J., Zhang, J.E. ~A unified intrinsic functional expansion theory for solitary waves.\\
\indent {\it Acta Mech Sinica}  21, 1-15 (2005).
\vskip 0.5mm
\noindent [9] von K$\acute{a}$rm$\acute{a}$n, Th. \& Burgers, J.M. General Aerodynamic Theory - Perfect fluids,  In {\it Aerodynamic\\ 
\indent Theory, Vol. II} (Ed. W.F. Durand) Calif. Inst. Tech. ~(1943).
\vskip 0.5mm
\noindent [10] Lighthill, M.J. A new approach to thin aerofoil theory.  {\it Aeron. Quart.} {\bf 3}, 193-210 ~(1951). 
\vskip 0.5mm
\noindent [11] Wu, Th. Y. \& Chwang, A.T. ~Double-body flow theory -- a new look at the classical problem.  In\\
\indent  {\it Tenth Symp. on Naval Hydrodynamics.} ONR 89-106 Dep't of the Navy, Washington DC (1974).
\vskip 0.5mm
\noindent [12] Wu, Th.Y. ~A nonlinear theory for unsteady flexible wing. {\it J. Eng. Math.} {\bf 58}, 279-287 (2007). 

\end{document}